\documentclass[a4paper,english,journal,onecolumn,draftcls,12pt]{IEEEtran}

\usepackage[T1]{fontenc}
\usepackage[latin9]{inputenc}
\usepackage{color}
\usepackage{amsthm, amsmath, amssymb}
\usepackage{esint}
\usepackage[numbers,sort&compress]{natbib}        
\usepackage{array}
\usepackage{datetime}
\usepackage{calc}
\usepackage{multicol,multirow}
\usepackage{bm,dsfont}
\usepackage{babel}
\theoremstyle{plain}
\newtheorem{theorem}{Theorem}
\newtheorem{lemma}[theorem]{Lemma}
\newtheorem{proposition}[theorem]{Proposition}

\newtheorem{corollary}[theorem]{Corollary}

\theoremstyle{definition}
\newtheorem{definition}[theorem]{Definition}
\newtheorem{remark}[theorem]{Remark}
\newtheorem*{remark*}{Remark}

\newcommand\otc[2]{%
\ifCLASSOPTIONonecolumn
\textcolor[rgb]{0.50,0.00,1.00}{}{#1}%
\else 
\textcolor[rgb]{0.00,0.59,0.00}{}{#2}%
\fi}

\usepackage{tikz,pgfplots}
\pgfplotsset{compat=newest}
\usetikzlibrary{external}

\newif\iffastfig

\newlength\figwidth
\newlength\figheight

\ifCLASSOPTIONonecolumn
    \usepackage{flcond-s}
\else
    \usepackage{flcond-d}
\fi

\begin{document}
\newlength{\norlen} \setlength{\norlen}{0.25ex} 

\global\long\def\st{\text{ s.t. }}
\global\long\def\OPmin#1#2#3{\begin{array}{cl}
 \min_{#1}  &  #2\\
\st &  #3 
\end{array}}

\global\long\def\norm#1{\left|\hspace{-\norlen}\left|\hspace{-\norlen}\left|#1\right|\hspace{-\norlen}\right|\hspace{-\norlen}\right|}
\global\long\def\vnorm#1{\left\Vert #1\right\Vert }
 \global\long\def\normi#1{\norm{#1}_{\infty}}

\global\long\def\ssum{{\textstyle \sum}}

\global\long\def\noise{\eta}
 \global\long\def\sinr{\tau}

\global\long\def\a{\mathbf{a}}
\global\long\def\b{\mathbf{b}}
\global\long\def\c{\mathbf{c}}
\global\long\def\d{\mathbf{d}}
\global\long\def\f{\mathbf{f}}
\global\long\def\fo{\mathbf{f}_{0}}
\global\long\def\g{\mathbf{g}}
\global\long\def\go{\mathbf{g}_{0}}
\global\long\def\h{\mathbf{h}}
\global\long\def\p{\mathbf{p}}
\global\long\def\q{\mathbf{q}}
\global\long\def\s{\mathbf{s}}
\global\long\def\v{\mathbf{v}}
\global\long\def\x{\mathbf{x}}
\global\long\def\xa{\mathbf{x^{\star}}}
\global\long\def\y{\mathbf{y}}

\global\long\def\qr{\mathfrak{q}}

\global\long\def\i{\mathcal{I}}
\global\long\def\ib{\bm{\i}}

\global\long\def\fx{\mathbf{f}\left(\x\right)}
\global\long\def\fox{\mathbf{f}_{0}\left(\x\right)}
 \global\long\def\ox{\left(\x\right)}
\global\long\def\pcc{\bm{\kappa}}

\global\long\def\I{\mathbf{I}}
 \global\long\def\0{\mathbf{0}}
\global\long\def\1{\mathds{1}}
 \global\long\def\R{\Re}
 \global\long\def\Rn{\R^{n}}

\global\long\def\S{\mathbf{S}}
\global\long\def\A{\mathbf{A}}
\global\long\def\B{\mathbf{B}}
\global\long\def\C{\mathbf{C}}

\global\long\def\iset{\mathcal{N}}
 \global\long\def\iseti{\mathcal{A}}
 \global\long\def\isete{\mathcal{B}}
 \global\long\def\D{\mathcal{D}}

\author{Martin~Jakobsson and Carlo~Fischione\thanks{
The authors are with the Department of Automatic Control, KTH Royal Institute of Technology, Sweden. Email: \{mjakobss,carlofi\}@kth.se.}
\thanks{The work of the authors is sponsored by the EU projects NoE HYCON2 and Hydrobionets.}}

\title{Optimality of Radio Power Control via \\
 Fast-Lipschitz Optimization}
\maketitle
\begin{abstract}
Fixed point algorithms play an important role to compute feasible solutions to the radio power control problems
in wireless networks. Although these algorithms are shown to converge
to the fixed points that give feasible problem solutions, the solutions often lack notion of problem optimality.
This paper reconsiders well known fixed point algorithms such as those with
standard and type-II standard interference functions, and investigates
the conditions under which they give optimal power control solutions by the recently proposed Fast-Lipschitz optimization
framework. When the qualifying conditions of Fast-Lipschitz optimization apply, it is established that the
fixed points are the optimal solutions of radio power optimization problems. The
analysis is performed by a logarithmic transformation of variables that gives problems treatable within the Fast-Lipschitz framework. It is shown how the logarithmic problem constraints are contractive by the standard or type-II standard assumptions on the
original power control problem, and how a set of cost functions fulfill the Fast-Lipschitz
qualifying conditions. The analysis on non monotonic interference function allows to establish a new qualifying condition for Fast-Lipschitz optimization. The results are illustrated
by considering power control problems with standard interference function,
problems with type-II standard interference functions, and a case
of sub-homogeneous power control problems. It is concluded that Fast-Lipschitz optimization
may play an important role in many resource allocation problems in
wireless networks.
\end{abstract}

\section{Introduction}
\label{sec:introduction}
Radio power control is one of the essential radio resource management techniques in wireless networks. The power control problem faces a tradeoff between saving power and having high enough level of power. It is important to control the transmit radio powers to avoid interferences to undesired receivers and save the energy of the transmitters. Meanwhile, it is important also to use adequate levels of power to make sure that the transmitted signals can overcome the attenuation of the wireless channel and the interference caused by other transmitters.

Power control in wireless communication is a particularly successful
instance of distributed optimization over networks.  Specifically,
in wireless networks a link is associated to one pair of nodes where
a node is a transmitter and the other is a receiver. Suppose there are $n$ transmitter-receiver pairs. Let $p_{i}$
be the radio power of transmit node $i$, for $i=1,\ldots,n$. Note that the index $i$ is used both for a transmitter and a receiver, so {\em transmitter $i$} and {\em receiver $i$} are two different nodes that are paired to communicate. One of the simplest
examples of radio power control considers the gain matrix $\mathbf{G}$,
where $G_{ij}$ is the channel attenuation from the transmit node
$j$ to the receiver node $i$. In addition to the useful signal $G_{ii}p_{i}$
from transmitter $i$, receiver $i$ will also receive a background
noise $\noise_{i}$ plus the interference of all other transmitters, $\sum_{j\neq i}G_{ij}p_{j}$.
For the communication attempt of transmit node $i$ to be successful, the signal
to (interference and) noise ratio (SNR) at receiver node $i$ must be higher than some threshold
$\sinr_{i}$,
\begin{equation}
\frac{G_{ii}p_{i}}{\sum_{j\neq i}G_{ij}p_{j}+\noise_{i}}\ge\sinr_{i}.\label{eq: SNR requirment}
\end{equation}
 If the transmit power of all links are collected in the vector \mbox{$\ensuremath{\p=[p_{1},p_{2},\dots,p_{n}]^{T},}$}
the requirement~\eqref{eq: SNR requirment} can be rewritten as
\begin{equation}
p_{i}\ge\i_{i}(\p) \triangleq \frac{\sinr_{i}}{G_{ii}}\left(\ssum_{j\neq i}G_{ij}p_{j}+\noise_{i}\right).\label{linear interference function}
\end{equation}
Using the vocabulary of \cite{Yates95}, we will refer to $\i_{i}(\p)$
as the \emph{interference function} of transmitter $i$. This affine
version of $\i_{i}(\p)$ is the simplest and best studied type of
interference function, and it is often the basis for extensions or
modifications by other types of interference functions. The focus on
achieving some minimum SNR in~\eqref{eq: SNR requirment} is justified
because many other measures of the quality of service are increasing
functions of the SNR \cite{Meshkati07}.

There are a number of ways of using the interference function in setting
power control optimization problems:
\begin{itemize}
\item maximization of the SNR (i.e., quality of service) of the
network, subject to power constraints;
\item minimization of the power consumption subject to SNR constraints;
\item maximization of some network utility function (e.g., throughput) of the
network, subject to power constraints.
\end{itemize}
Early works on distributed power control in wireless networks have followed
the first approach, and try to maximize the smallest SNR of the
network \cite{Zander,Grandhi94}. With the inclusion of receiver noise
in \cite{Foschini1993}, the focus has shifted to the second approach,
with the goal of minimizing the radio powers $p_{i}$ while maintaining
a minimum SNR $\sinr_{i}$ at each receiver, i.e.,

\begin{equation} \label{eq power problem}
\begin{array}{l l}
\min_{\p} & \quad\p  \\
{\rm s.t.} & \quad p_{i}\ge\i_{i}(\p)\quad\forall i.
\end{array}
\end{equation}

This line of work has later been generalized to the framework
of \emph{standard} interference functions by Yates~\cite{Yates95}.
When problem~\eqref{eq power problem} above is feasible, and the
functions $\i_{i}(\p)$ are standard (see Definition~\ref{def: standard}),
the unique optimal solution to \eqref{eq power problem} is given by the fixed point of the iteration
\begin{equation}
p_{i}^{k+1}:=\i_{i}(\p^{k}),\label{iterations}
\end{equation}
 or, in vector form, $\p^{k+1}:=\ib(\p^{k})$ where
\[
\ib(\p^{k})\triangleq\begin{bmatrix}\i_{1}(\p^{k}) & \dots & \i_{n}(\p^{k})\end{bmatrix}^{T}.
\]
The computation of the optimal solution for problem~\eqref{eq power problem}
by these iterations is much simpler than using the classical parallelization
and decomposition methods of distributed optimization \cite{Bertsekas1997}.
This is because there is no longer a need to centrally collect, compute
and redistribute the primal and dual variables of the problem due to that
$\i_{i}(\p^{k})$ can be known or estimated locally at receiver $i$~\cite{Foschini1993,Yates95}.
Even in a centralized setting, iteration \eqref{iterations} is simpler
than traditional distributed optimization methods, since no dual variables
need to stored and manipulated. The iterations require only that every receiver
node successively updates using local knowledge of the function (interference function) of other nodes' current
decision variables (radio powers). Another advantage is that convergence
is obtained even though such a knowledge is delayed, i.e., the decision
variables $p_{j}^{k}$ of other nodes are updated with some delay~\cite{Bertsekas1997,Mitra93}.

Give the advantages mentioned above, there is a number of studies in the literature where radio power
control algorithms have been proposed by considering iterations similar
to~\eqref{iterations} \cite{Leung2004,Sung99,Sung2005,Ulukus98,Luo05,Biguesh11}.
In these approaches, the interference functions are not necessarily standard,
and the focus has been in studying the convergence of iteration~\eqref{iterations} to a fixed point
rather than the meaning of the fixed point in terms of optimality for problem~\eqref{eq power problem}. One of the reasons is that optimal power control problems of the form~\eqref{eq power problem} with or without standard interference functions may be non convex and non linear, which makes it very hard the characterization of optimality. Therefore, the natural question that we would like to answer in the paper is
whether there are conditions ensuring the optimality of existing fixed
point radio power control algorithms.

 Fast-Lipschitz optimization is a recently proposed framework that
is motivated by such a question~\cite{Fischione11,Jakobsson13}.
In particular, Fast-Lipschitz optimization is a natural generalization
of the interference function approach on how to solve distributed
optimization problems over wireless networks by using fixed point
iterations similar to \eqref{iterations}. Fast-Lipschitz does not
need constraints that are standard, and can have a more general objective
functions than the one in problem~\eqref{eq power problem}. The
main characteristic of a Fast-Lipschitz problem is that the optimal
point is given by the fixed point of the constraints, a result that is in general very difficult to establish. In this paper
we investigate under which conditions a general power control problem
falls under the Fast-Lipschitz framework, whereby the fixed points
of iteration~\eqref{iterations} are optimal.

The reminder of this paper is organized as follows: In Subsection
\ref{sub: Related-work} we discuss the related work and in Subsection
\ref{sub: Notation} notation is introduced. The problem formulation
is given in Section~\ref{sec:Problem-Formulation}, and for the sake
of self-containment we give a brief definition of Fast-Lipschitz optimization
in Section~\ref{sec:Fast-Lipschitz-Optimization}. In Section~\ref{sec: TSS preliminary}
we present preliminary results on two-sided scalable functions, which
we use to examine the relation between standard and type-II standard
functions in sections \ref{sec: standard is FL} and \ref{sec: type-II is FL}
respectively. Section \ref{sec:Absolutely-subhomogeneous-interf}
gives an example of Fast-Lipschitz optimization applied to a problem
that is neither standard, nor type-II standard. Finally, the paper
is concluded in Section~\ref{sec:Conclusions-and-future}.

\subsection{Related work \label{sub: Related-work}}

The iterative methods to solve radio power control problems are a special case of parallel and distributed
computation theory. There is a rich line of research on
distributed iterative methods, with the corner stones \cite{Tsitsikilis,Bertsekas1997}
by Tsitsiklis and Bertsekas. Most of the recent work focuses on convex
problems, where duality and decomposition techniques can be used to
distribute the computations over the involved nodes or agents of the
network (see, e.g., \cite{Palomar2007} for a discussion of different
methods). A framework that recently has attracted substantial attention
is the Alternating method of Lagrangian Multipliers (ADMM) \cite{Boyd2011}, which has been particularly successful for
optimization problems in learning theory with huge data sets. These
methods require a central entity that coordinates the nodes and their
optimization subproblems. The problems are distributed in a computational
sense, meaning that every node makes the computations coordinate by some central entity, rather than decentralized from an organizational point of
view.

Decentralized solution methods are addressed by consensus methods,
where all nodes are peers and compute the solution of network optimization
problems by exchanging information with their local neighbors (see,
e.g., \cite{Boyd2006,Olshevsky2009,Nedic2010}). In \cite{Nedic2010},
each node has a local cost function and a local constraint set, both of which are assumed
convex. The local problems are coupled through a common variable,
and the global objective is to minimize the sum of all local costs.
These powerful optimization methods are not easy to apply to the optimization
of wireless networks, due to the slowness of the convergence of the
message passing procedure. For example, consensus methods may converge
with some hundreds of message exchanges that
would take more time than the time needed to compute the optimal transmit power compared to the coherence time of the wireless channel.
This means that when the solution would be computed by consensus methods, it is outdated.

In power control problems we often have in wireless networks, the global cost is not necessarily separable,
nor convex, and each node controls only a subset of the variables.
This is coherent with the lines of Yates' framework \cite{Yates95} and
the algorithms that are standard (e.g., \cite{Zander,Mitra93,Foschini1993,Grandhi94}).
In fact, this paper will show that standard algorithms are encompassed in
the Fast-Lipschitz framework. In \cite{Leung2004}, Yates' framework
is generalized to cover also some discrete implementations (e.g.,
\cite{Sung99}, where the power updates are of a fixed size). Extensions of Yates' framework have been proposed by Sung and Leung
\cite{Sung2005}. They consider opportunistic algorithms that are
not standard by Yates' definition. Instead, they introduce \emph{type-II
standard} functions and the more general \emph{two-sided scalable}
interference functions. These functions are shown to have the same
fixed point properties as Yates' standard functions, i.e., problem~\eqref{eq power problem}
with type-II standard constraints can be solved through repeated iterations
of the constraints.

Further extension of the interference function
framework have been proposed in~\cite{Boche08,Hamid}. Specifically,\cite{Boche08}
considers standard functions with a small modification, where the \emph{scalability} property (see section~\eqref{sec: standard is FL})
is replaced by the \emph{scale invariance}, i.e., \emph{$\i(c\p)=c\i(\p)$.}. It is shown how these functions can be interpreted as level curves
of closed comprehensive sets. \cite{Hamid} replaces the scalability
assumption of the standard framework by a requirement of weighted
max-norm contraction. This allows to derive statements on the fixed
points and convergence speed of iteration~\eqref{iterations}. However,
neither of these extensions are concerned with notions of optimality
of problems in the form of \eqref{eq power problem}.

The algorithms above assume perfect knowledge of the
interference functions. In \cite{Ulukus98} the convergence when some
of the measurements required to evaluate the interference functions
(e.g., the SNR samples) are stochastic and noisy or inaccurate is
studied. \cite{Luo05} shows the convergence of a general class of
stochastic power control algorithms, given a suitable set of power
update damping step lengths, and choice of these lengths are decentralized
and improved in \cite{Biguesh11}.

A different approach to power control is based on game theory, e.g.,
\cite{Hongbin98,Saraydar02,Xiao03,Meshkati07,Feng13}. \cite{Hongbin98}
introduced an economic framework and modeled the power control problem
as a non-cooperative game, where each mobile (node
of the network) selfishly tries to maximize its local utility. The
corresponding power updates \eqref{iterations} are then the best
response of each mobile, and the fixed point corresponds to a Nash
equilibrium. \cite{Saraydar02} investigates how this Nash equilibrium
is affected by pricing of transmit powers. The game theoretic framework
is flexible enough to model also cognitive radio networks with primary
and secondary users \cite{Feng13}. An open problem in this line of
research is that the resulting Nash equilibria typically do not correspond
to a social optimum, meaning that there are other power configurations
where all users are better off, or some global utility (such as total
throughput) is higher.

In the power control algorithms mentioned above, the focus is not
about the optimally of radio power control, but the convergence and
existence of equilibria in distributed power updatings. In this paper, we consider the recently proposed theory
of Fast-Lipschitz optimization to establish the optimally
of power control algorithms. In particular, Fast-Lipschitz optimization
is related to other techniques that replace the most common assumption
of convexity with other conditions. Examples include the framework
of monotonic optimization (e.g., \cite{Tuy2000,Bjornson2012,Zhang13}),
where monotonic properties of the objective and the constraints allow
for efficient solutions to the problem, or the framework of abstract
optimization \cite{Notarstefano2011,burger12} that generalizes linear
programming in the sense that problems are solved by determining the
subsets of the constraints that define the solution.

\subsection{Notation \label{sub: Notation}}

Vectors and matrices are denoted by bold lower and upper case letters,
respectively. The components of a vector $\x$ are denoted $x_{i}$
or $[\x]_{i}$. Similarly, the elements of the matrix $\mathbf{A}$
are denoted $A_{ij}$ or $\left[\mathbf{A}\right]_{ij}$ . The transpose
of a vector or matrix is denoted $\cdot^{T}$. $\I$ and $\1$ denote
the identity matrix and the vector of all ones. A vector or matrix
where all elements are zero is denoted by $\0$.

The gradient of a function is defined as $\left[\nabla\fx\right]_{ij}=\partial f_{j}\ox/\partial x_{i}$,
whereas $\nabla_{i}\fx$ denotes the $i$th row of $\nabla\fx$. Note
that $\nabla\fx^{k}=\left(\nabla\fx\right)^{k}$, which has not to
be confused with the $k$th derivative. The spectral radius is denoted
$\rho(\cdot)$. Vector norms are denoted $\vnorm{\cdot}$ and matrix
norms are denoted $\norm{\cdot}$. Unless specified $\vnorm{\cdot}$
and $\norm{\cdot}$ denote arbitrary norms. $\normi{\mathbf{A}}=\max_{i}\sum_{j}|A_{ij}|$
is the norm induced by the $\ell_{\infty}$ vector norm, where $\vnorm{\x}_{\ell_{\infty}}=\max_{i}\left|x_{i}\right|$.
These matrix norm definitions are coherent with \cite{Horn85}.

All inequalities in this paper are intended \emph{element-wise}, i.e.,
$\A\ge\B$ means $A_{ij}\ge B_{ij}$ for all $i,j$. We will also
use the element-wise natural logarithm $\ln\x=\left[\ln x_{1},\dots\,,\ln x_{n}\right]^{T}$
and the element-wise exponential $e^{\x}=\left[e^{x_{1}},\dots\,,e^{x_{n}}\right]^{T}$

\begin{remark*} The notation $\i$ for interference functions does
not follow the notational assumptions above. However, we will keep
the notation to harmonize with existing literature~\cite{Grandhi94,Foschini1993,Yates95,Boche08}.
\end{remark*}

\section{Problem Formulation \label{sec:Problem-Formulation}}

We investigate a general form of the power minimization problem mentioned
in the introduction section and having the general form
\begin{equation}\label{op: power control}
\OPmin{\p}{\pcc\left(\p\right)}{\p\ge\ib\left(\p\right).}
\end{equation}
Throughout the paper we assume that $\pcc(\p)$ and $\ib(\p)$ are
differentiable. The cost function $\pcc\left(\p\right)$ can be scalar
or vector valued. Examples are $\pcc\left(\p\right)=\p$ or $\pcc\left(\p\right)=\p^{T}\1$.
In practice, the powers must be positive and there is a maximum power
that each transmitter can generate. Therefore, we will implicitly
assume that there are the natural constraints $\p\in\D_{\p}=\left\{ \p\,:\,\p_{\min}\le\p\le\p_{\max}\right\} $,
where $\p_{\min}\ge\0$ and $\p_{\max}$ are given constants.

The main problem this paper is concerned with, is when the iterations
\begin{equation}
\p^{k+1}:=\ib\left(\p^{k}\right)\label{eq: fixed point iterations}
\end{equation}
solve optimization problem~\eqref{op: power control}. We show under
which conditions the general power control problem (5) is Fast-Lipschitz,
which will allow us to establish the optimality of iterations (6).
In particular, when problem \eqref{op: power control} is Fast-Lipschitz,
then the iterations \eqref{eq: fixed point iterations} will converge
to $\p^{\star}=\ib(\p^{\star})$ and $\p^{\star}$ is optimal for
problem \eqref{op: power control}.
\begin{remark}
Note that the formulation of the iterations \eqref{eq: fixed point iterations}
is synchronous, i.e., every node must finish the computations and
communications of round $k$ before the next round $k+1$ can start.
The algorithm we consider also converges asynchronously, under the
assumption of bounded delays \cite{Mitra93,Sung2005,Bertsekas1997,Fischione11}.
However, since convergence properties are not the main focus of this
paper, we restrict ourselves to the less cumbersome synchronous notation
of \eqref{eq: fixed point iterations}.
\end{remark}

\section{Fast-Lipschitz Optimization \label{sec:Fast-Lipschitz-Optimization}}

For the sake of self containment and the need for introducing a preliminary
result, we now give a brief formal definition of Fast-Lipschitz problems.
For a thorough discussion of Fast-Lipschitz properties we refer the
reader to \cite{Fischione11,Jakobsson13}.
\begin{definition}
\label{def: FL form}A problem is said to be on\emph{ Fast-Lipschitz
form} if it can be written as
\begin{equation}
\begin{array}{cl}
\max & \fo(\x)\\
\text{s.t.} & x_{i}\le f_{i}(\x)\quad\forall i\in\iseti\\
 & x_{i}=f_{i}(\x)\quad\forall i\in\isete,
\end{array}\label{eq FLform}
\end{equation}
 where
\end{definition}
\setlength{\itemsep}{1ex}
\begin{itemize}
\item $\fo:\Rn\rightarrow\Re^{m}$ is a differentiable scalar ($m=1$) or
vector valued ($m\ge2$) function.
\item $\iseti$ and $\isete$ are complementary subsets of $\{1,\dots,n\}$.
\item For all $i$, $f_{i}\,:\,\Rn\to\Re$ is a differentiable function.
\end{itemize}
 From the individual constraint functions we form the vector valued
function $\f:\Rn\rightarrow\Rn$ as \ensuremath{\f(\x)=\begin{bmatrix}f_{1}(\x) & \cdots & f_{n}(\x)\end{bmatrix}^{T}.}

\begin{remark}
For the rest of the paper, we will restrict our attention to a \emph{bounding
box} \ensuremath{\D=\left\{ \x\in\Rn\,|\,\a\le\x\le\b\right\} .}
 We assume $\D$ contains all candidates for optimality and that $\f$
maps $\D$ into $\D$, \ensuremath{\f:\D\to\D.}
 This box arises naturally in practice, since any real-world
decision variable must be bounded.
\end{remark}

\begin{definition}
A problem is said to be \emph{Fast-Lipschitz} when it can be written
on Fast-Lipschitz form and admits a unique Pareto optimal solution
$\xa$, defined as the unique solution to the system of equations
\[
\xa=\f(\xa).
\]
\end{definition}

A problem written on Fast-Lipschitz form is not automatically Fast-Lipschitz.
The appendix provides qualifying conditions which, when fulfilled,
guarantee that problem~\eqref{eq FLform} is Fast-Lipschitz (see
Table~\ref{tab:oldcond} and Theorem~\ref{theorem: Main theorem - New}).

The framework of Fast-Lipschitz optimization is formulated for maximization
problems. Through a change of variables, any minimization problem
can be formulated as an equivalent maximization problem. This is useful
when dealing with power control problems, which are normally written
by minimization. The following lemma shows how minimization problems
fit into the Fast-Lipschitz framework.
\begin{lemma}
[Fast-Lipschitz minimization] Consider
\begin{equation}
\begin{array}{cl}
\min & \g_{0}(\x)\\
\text{s.t.} & x_{i}\ge g_{i}(\x)\quad\forall i\in\iseti\\
 & x_{i}=g_{i}(\x)\quad\forall i\in\isete,
\end{array}\label{eq FLminForm}
\end{equation}
 where $\g_{0}(\x)$, $\g(\x)=[g_{i}(\x)]$, $\iseti$ and $\isete$
fulfill the assumptions of Definition~\ref{def: FL form}. Then,
problem \eqref{eq FLminForm} is Fast-Lipschitz if $\g_{0}(\x)$ and
$\g(\x)$ fulfill the qualifying conditions. \end{lemma}
\begin{IEEEproof}
Let $\x=-\y$ and form the equivalent problem
\begin{equation}
\begin{array}{cl}
\max & -\g_{0}(-\y)=\f_{0}(\y)\\
\text{s.t.} & y_{i}\le-g_{i}(-\y)=f_{i}(\y)\quad\forall i\in\iseti\\
 & y_{i}=-g_{i}(-\y)=f_{i}(\y)\quad\forall i\in\isete.
\end{array}\label{eq FLminForm2}
\end{equation}
 In order to check the qualifying conditions one needs $\nabla\f_{0}(\y)$
and $\nabla\f(\y)$. But
\otc{%
\begin{align*}
(\nabla_{\y}\f_{0}(\y))^{T}
&
=\frac{\partial\f_{0}(\y)}{\partial\y}
=\frac{\partial\left(-\g_{0}(\x)\right)}{\partial\x}\frac{\partial\x}{\partial\y}
=-\frac{\partial\g_{0}(\x)}{\partial\x}(-1)
=(\nabla_{\x}\g_{0}(\x))^{T},
\end{align*}
}{%
\begin{align*}
(\nabla_{\y}\f_{0}(\y))^{T}
&
=\frac{\partial\f_{0}(\y)}{\partial\y}
=\frac{\partial\left(-\g_{0}(\x)\right)}{\partial\x}\frac{\partial\x}{\partial\y}
=-\frac{\partial\g_{0}(\x)}{\partial\x}(-1)
\\&
=(\nabla_{\x}\g_{0}(\x))^{T},
\end{align*}
}
 and analogously, $\nabla_{\y}\f(\y)=\nabla_{\x}\g(\x)$. Since $\g_{0}(\x)$
and $\g(\x)$ fulfill the qualifying conditions, the equivalent problem
\eqref{eq FLminForm2} is Fast-Lipschitz.
\end{IEEEproof}
We are now in the position of introducing the core contribution of
the paper in the following section.

\section{Two-sided scalable problems and Fast-Lipschitz optimization}

In this section we examine the relations between standard functions
(in \ref{sec: standard is FL}), type-II standard functions (in \ref{sec: type-II is FL}),
and Fast-Lipschitz optimization. This will allow us to establish the
core results that: 1) all standard power control problems are Fast
Lipschitz; 2) there exist non standard interference functions whose
fixed point is the optimal of a power control problem; 3) type II
power control iterations are optimal within some conditions. We begin
by briefly recalling the concept of two-sided scalability, which we
will use to put standard and type-II standard functions in the Fast-Lipschitz
framework.

\subsection{Preliminary results on two-sided scalable functions}

\label{sec: TSS preliminary}

This subsection presents preliminary results for the upcoming sections
on standard and type-II standard problems. The main result of this
section is Lemma~\ref{lem: norm_1 <1}, which will allow us to establish
the contractivity of standard and type-II standard functions. Contractivity is one of the main ingredients for Fast-Lipschitz optimization.
\begin{definition}
[\cite{Sung2005}]\label{def: two-sided scalable}A function $\ib\left(\p\right)$
is \emph{two-sided scalable }if for all $c>1$ and all $(1/c)\p\le\q\le c\p,$
it holds that
\begin{equation}
(1/c)\ib(\p)<\ib(\q)<c\ib(\p).\label{eq: two-sided scalable}
\end{equation}
\end{definition}
\begin{proposition}
[{\cite[Prop. 4]{Sung2005}}]\label{prop: standard is 2sided-scalable}
If a function is either standard or type-II standard, then it is also
two-sided scalable.
\end{proposition}

The following lemma shows that two-sided scalable functions are shrinking maps with gradients
of one-norm less than one. The lemma  is based on the logarithmic transformations proposed
in \cite{Moller2013}. It will be used in the main results of sections
\ref{sec: standard is FL} and \ref{sec: type-II is FL}.
\begin{lemma}
\label{lem: norm_1 <1}Let $\x\triangleq\ln\p$ and $\f(\x)\triangleq\ln\ib(\x)$.
If $\ib\left(\p\right)$ is two-sided scalable, then
\[
\vnorm{\f\left(\x\right)-\f\left(\y\right)}_{\infty}<\vnorm{\x-\y}_{\infty}
\]
for all $\x,\y$, and
\[
\norm{\nabla\fx}_{1}<1.
\]
\end{lemma}
\begin{IEEEproof}
By \cite[Lemma 7]{Sung2005}, two-sided scalability implies
\[
\max_{i}\left\{ \max\left\{ \frac{\i_{i}\left(\x\right)}{\i_{i}\left(\y\right)},\frac{\i_{i}\left(\y\right)}{\i_{i}\left(\x\right)}\right\} \right\} <\max_{i}\left\{ \max\left\{ \frac{x_{i}}{y_{i}},\frac{y_{i}}{x_{i}}\right\} \right\}
\]
for all $i$. Since the logarithm is strictly increasing, this is
equivalent to
\otc{%
\[
\max_{i}\left\{ \max\left\{ \ln\left(\frac{\i_{i}\left(\p\right)}{\i_{i}\left(\q\right)}\right),\ln\left(\frac{\i_{i}\left(\q\right)}{\i_{i}\left(\p\right)}\right)\right\} \right\}
<\max_{i}\left\{ \max\left\{ \ln\left(\frac{p_{i}}{q_{i}}\right),\ln\left(\frac{q_{i}}{p_{i}}\right)\right\} \right\} .
\]
}{%
\begin{multline*}
\max_{i}\left\{ \max\left\{ \ln\left(\frac{\i_{i}\left(\p\right)}{\i_{i}\left(\q\right)}\right),\ln\left(\frac{\i_{i}\left(\q\right)}{\i_{i}\left(\p\right)}\right)\right\} \right\} \\
<\max_{i}\left\{ \max\left\{ \ln\left(\frac{p_{i}}{q_{i}}\right),\ln\left(\frac{q_{i}}{p_{i}}\right)\right\} \right\} .
\end{multline*}
}
Inserting $\p=e^{\x}$ and $\q=e^{\y}$ gives
\[
\max_{i}\left\{ \left|\ln\i_{i}\left(e^{\x}\right)-\ln\i_{i}\left(e^{\y}\right)\right|\right\} <\max_{i}\left\{ \left|\ln e^{x_{i}}-\ln e^{y_{i}}\right|\right\} ,
\]
\[
\Leftrightarrow\left|f_{i}\left(\x\right)-f_{i}\left(\y\right)\right|<\left|x_{i}-y_{i}\right|.
\]
Since this holds for all components $i$, we have
\[
\vnorm{\fx-\f\left(\y\right)}_{\infty}<\vnorm{\x-\y}_{\infty}
\]
by definition of norm infinity.

Denote $\v=\arg\max_{\vnorm{\bf u}_{\infty}=1}\vnorm{\nabla\fx^{T}\bf u}_{\infty}$.
By definition, we have $\vnorm{\v}_{\infty}=1$ and
\[
\vnorm{\nabla\fx^{T}\v}_{\infty}=\norm{\nabla\fx^{T}}_{\infty}=\norm{\nabla\fx}_{1}.
\]
By defining $\y=\x+\epsilon\v$, with $\epsilon$ positive scalar, we have
 \otc{%
\begin{align*}
1 & >\frac{\vnorm{\f\left(\y\right)-\fx}_{\infty}}{\vnorm{\y-\x}_{\infty}}
=\frac{\vnorm{\f\left(\x+\epsilon\v\right)-\fx}_{\infty}}{\vnorm{\epsilon\v}_{\infty}}
=\vnorm{\frac{\f\left(\x+\epsilon\v\right)-\fx}{\epsilon}}_{\infty},
\end{align*}
}{%
\begin{align*}
1 & >\frac{\vnorm{\f\left(\y\right)-\fx}_{\infty}}{\vnorm{\y-\x}_{\infty}}=\frac{\vnorm{\f\left(\x+\epsilon\v\right)-\fx}_{\infty}}{\vnorm{\epsilon\v}_{\infty}}
\\& =
\vnorm{\frac{\f\left(\x+\epsilon\v\right)-\fx}{\epsilon}}_{\infty},
\end{align*}}
and in the limit $\epsilon\to0$,
\otc{%
\begin{align*}
1
&
> \lim_{\epsilon\to0}\frac{\vnorm{\f\left(\y\right)-\fx}_{\infty}}{\vnorm{\y-\x}_{\infty}}
=\lim_{\epsilon\to0}\vnorm{\frac{\f\left(\x+\epsilon\v\right)-\fx}{\epsilon}}_{\infty}
=\vnorm{\nabla\fx^{T}\v}_{\infty}=\norm{\nabla\fx}_{1}.
\end{align*}}{%
\begin{align*}
1
&
> \lim_{\epsilon\to0}\frac{\vnorm{\f\left(\y\right)-\fx}_{\infty}}{\vnorm{\y-\x}_{\infty}}
=\lim_{\epsilon\to0}\vnorm{\frac{\f\left(\x+\epsilon\v\right)-\fx}{\epsilon}}_{\infty}
\\&
=\vnorm{\nabla\fx^{T}\v}_{\infty}=\norm{\nabla\fx}_{1}.
\end{align*}}
This concludes the proof.
\end{IEEEproof}
The result shows that two-sided scalable functions are shrinking maps,
but it does not establish the amount of slack in $\norm{\nabla\fx}_{1}<1$.
This slack is useful for several reasons. The first
reason is that any amount of slack makes the function $\fx$
a contraction, thereby guaranteeing a unique fixed point. This is
assumed in the Fast-Lipschitz qualifying conditions, e.g., in \qcc{Qk.cont}.
The lack of knowledge on the amount of slack is not a problem in practice, since the bounds $\p_{\min}\le\p\le\p_{\max}$
form a closed bounded region of $\Rn$. Therefore, $\norm{\nabla\fx}_{1}$
attains a minimum (call this value $\alpha$) in that region, so $\norm{\nabla\fx}_{1}\le\alpha<1$
for those $\x$ that are of interest. Secondly (and more importantly),
the qualifying conditions other than \qc{Q1} typically require $\norm{\nabla\fx}_{\infty}<c$
for some $c\in(0,1]$. Lemma~\ref{lem: norm_1 <1} is therefore only
of use if $c=1$. This special case is exploited in Section~\ref{sec: type-II is FL}. Now, based on this lemma, we are in the position to show that power control problems~\eqref{op: power control} with standard interference functions are a special case and Fast-Lipschitz optimization.

\subsection{Standard functions and Fast-Lipschitz optimization}

\label{sec: standard is FL}

In this section we recall Yates' framework of standard functions and
show that a problem \eqref{op: power control}  with standard interference function constraints has an equivalent problem formulation
that is Fast-Lipschitz.
\begin{definition}
[\cite{Yates95}]\label{def: standard}The function $\ib(\p)$ is
\emph{standard} if for all $\p,\q\ge\0$, the following properties
are satisfied. \begin{subequations}
\begin{alignat}{2}
 & \textit{Positivity: } &  & \ib(\p)>\0\\
 & \textit{Monotonicity: } &  & \p\ge\q\Rightarrow\ib(\p)\ge\ib(\q)\label{eq: standard a}\\
 & \textit{Scalability: } &  & c>1\Rightarrow\ib(c\p)<c\ib(\p)\label{eq: standard b}
\end{alignat}
\end{subequations} The monotonicity requirement \eqref{eq: standard a}
can equivalently be formulated $\nabla\ib\left(\p\right)\ge\0\text{ for all }\p\ge\0.$

There is a relation between standard functions and two-sided scalable functions. Note that \eqref{eq: two-sided scalable} multiplied
by a positive scalar $c>0$ implies $\left(c^{2}-1\right)\ib\left(\p\right)>\0$,
so any two-sided scalable function must also be positive, i.e., $\ib\left(\p\right)>\0$
\cite[Lemma 6]{Sung2005}. If $\q$ in Definition~\ref{def: two-sided scalable}
is chosen as $\q=c\p$, inequalities \eqref{eq: two-sided scalable}
become
\begin{equation}
(1/c)\ib(\p)<\ib(c\p)<c\ib(\p),\label{eq: 2sided scalable}
\end{equation}
so a two-sided scalable function is always scalable \eqref{eq: standard b}.
Any two-sided scalable function $\ib\left(\p\right)$ is therefore
standard if $\nabla\ib\left(\p\right)\ge\0$.
\end{definition}

\begin{proposition}
[\cite{Yates95}] \label{thm: Standard function}Assume that the
power optimization problem~\eqref{eq power problem} is feasible.
Then, the standard interference function $\ib(\p)$ has a unique fixed
point $\p^{\star}$ that is the solution to~\eqref{eq power problem}.
\end{proposition}
{}

We now show that problem \eqref{eq power problem} with standard
interference constraints fall under the Fast-Lipschitz framework.
To this end we consider problem \eqref{op: power control}, with the
general cost function $\pcc\left(\p\right)$, which we assume be
differentiable. All Fast-Lipschitz qualifying conditions require that
the norm of the constraint function gradient be small enough. We will
show this through Lemma~\ref{lem: norm_1 <1}, wherefore we investigate
problem~\eqref{op: power control} after a change of variables.

To this end, we let $\x\triangleq\ln\p$ as
the logarithm of the power variables. This gives $\p\triangleq e^{\x}$
and the equivalent problems
\[
\OPmin{\x}{\pcc\left(e^{\x}\right)}{e^{\x}=\ib(e^{\x})}
\]
and, because the logarithm is strictly increasing,
\begin{equation}
\OPmin{\x}{\fox \triangleq \pcc\left(e^{\x}\right)}{\x\ge\fx\triangleq\ln\ib(e^{\x}).}\label{op: type1 fl}
\end{equation}
If problem~\eqref{op: type1 fl} is Fast-Lipschitz, then $\xa=\f\left(\xa\right)$
is the unique Pareto optimal point for \eqref{op: type1 fl}, whereby
\[
e^{\xa}=\p^{\star}=\ib\left(\p^{\star}\right)
\]
 is optimal for problem \eqref{op: power control}.

The following result shows how power control problems with standard constraint functions, if differentiable,
have an equivalent Fast-Lipschitz problem formulation.
\begin{theorem}
\label{thm: type-1 is FL-1}Consider problem~\eqref{op: power control}
and let $\ib\left(\p\right)$ be differentiable and standard. If
$\nabla\pcc(\p)\ge\0$ with non-zero rows, then the
equivalent problem~\eqref{op: type1 fl} is Fast-Lipschitz and
$\p^{\star}=\ib\left(\p^{\star}\right)$ is optimal in problem~\eqref{op: power control}.
\end{theorem}
\begin{IEEEproof}
\global\long\def\inlinenorm#1{|\hspace{-\norlen}|\hspace{-\norlen}|#1|\hspace{-\norlen}|\hspace{-\norlen}|}
We show that problem~\eqref{op: type1 fl} is Fast-Lipschitz by qualifying
condition \qc{Q1}. The gradients of problem~\eqref{op: type1 fl}
are given by \begin{subequations}\label{eq: gradients}
\begin{align}
 & \nabla\fox=\text{diag}\left(\p\right)\nabla\pcc\left(\p\right),\\
 & \nabla\fx=\text{diag}\left(\p\right)\nabla\ib(\p)\,\text{diag}\left(1/\ib\left(\p\right)\right),
\end{align}
 \end{subequations}where $\p=e^{\x}$ and $\left[\text{diag}\left(1/\ib\left(\p\right)\right)\right]_{ii}=1/\i_{i}\left(\p\right)$.
Since $\p=e^{\x}\ge\0$, $\ib(\p)>\0$ and $\nabla\ib(\p)\ge\0$,
the gradients~\eqref{eq: gradients} fulfill $\nabla\fx\ge\0$ and
$\nabla\fox\ge\0$ with non-zero rows (these are conditions \qcc{Q1.pos}
and \qcc{Q1.f0} respectively). If $\ib\left(\p\right)$ is standard,
it is also two-sided scalable by Proposition~\ref{prop: standard is 2sided-scalable},
so $\norm{\nabla\f\left(\x\right)}_{1}<1$ by Lemma~\ref{lem: norm_1 <1}
(condition \qcc{Q1.cont}). Problem~\eqref{op: type1 fl} is therefore
Fast-Lipschitz by qualifying condition~\qc{Q1}, and
$\x^{\star}=\f\left(\x^{\star}\right)$. Taking the exponential of
the previous relation gives $\p^{\star}=\ib\left(\p^{\star}\right)$.
This concludes the proof.
\end{IEEEproof}
While Proposition~\ref{thm: Standard function} states that the fixed
point of standard constraints minimize the powers in a Pareto sense,
i.e., $\pcc\left(\p\right)=\p$, Theorem~\ref{thm: type-1 is FL-1}
accepts any non-decreasing $\pcc\left(\p\right)$. The requirement
that $\nabla\pcc\left(\p\right)$ have non-zero rows simply means that
each variable $p_{i}$ has an effect on at least one component of
the cost at each $\p$. For scalar values cost functions, this is
the same as requiring $\pcc$ to be strictly increasing, $\nabla\pcc\left(\p\right)>\0.$

Theorem~\ref{thm: type-1 is FL-1} is not a generalization of Theorem~\ref{thm: Standard function}
in practice, since a minimization of $\p$ is equivalent to a minimization
of an increasing function of $\p$. The novelty here is instead that
standard problems falls within the broader class of Fast-Lipschitz
problems. Therefore, we can have non standard interference functions
that can lead to optimally by distributed iterative power control
algorithms.  In the next subsection we will continue to show how type-II
standard functions relates to Fast-Lipschitz optimization.

\subsection{Type-II standard functions and Fast-Lipschitz optimization}

\label{sec: type-II is FL}

As the standard functions are monotonically increasing, transmit nodes following
\eqref{eq: fixed point iterations} will always increase their power
when their transmission environment is worsened by higher interference. A receiver node experiencing a deep
fade will therefore need a very high transmit power, thereby increasing interference
for the other receiver nodes in the network. This is not a good strategy, for example, in
delay tolerant applications, where transmit nodes can adjust their transmission
rates and higher throughput can be achieved by prioritizing receiver nodes experiencing
low interference. One such strategy is to keep the signal-to-interference
product constant, which results in update functions \eqref{eq: fixed point iterations}
that are monotonically decreasing, and therefore not standard. This
is addressed in \cite{Sung2005}, where Sung and Leung extends Yates'
framework with type-II standard functions.
\begin{definition}
[\cite{Sung2005}]The function $\ib(\p)$ is \emph{type-II standard}
if for all $\p,\q\ge\0$, the following properties are satisfied:
\begin{subequations}\label{def: type2standard}
\begin{alignat}{2}
 & \textit{Type-II Monotonicity: } &  & \p\le\q\Rightarrow\ib(\p)\ge\ib(\q)\label{eq: type-II a}\\
 & \textit{Type-II Scalability: } &  & c>1\Rightarrow\ib(c\p)>(1/c)\ib(\p)\label{eq: t2 scalable}
\end{alignat}
\end{subequations}
\end{definition}
As in the case of standard functions, the monotonicity property \eqref{eq: type-II a}
can be written $\nabla\ib\left(\p\right)\le\0$ for all $\p$. Note
also from \eqref{eq: 2sided scalable} that all two-sided
scalable functions are type-II scalable \eqref{eq: t2 scalable} so
any two-sided scalable function $\i\left(\p\right)$ where $\nabla\i\left(\p\right)\le\0$
is also type-II standard. Type-II standard functions converge in
the same way as standard functions - if a fixed point $\p^{\star}$
exists, then iteration~\eqref{eq: fixed point iterations} converges
to $\p^{\star}$ \cite[Thm. 3]{Sung2005}.

When considering opportunistic algorithms, $\i_{i}\left(\p\right)$
no longer has the interpretation of {}``interference that receiver node
$i$ must overpower asking transmit node $i$ to use a power $p_{i}$ high enough''. There are no longer any explicit
constraints $\p\ge\ib\left(\p\right)$ underlying the algorithm, and
$\ib$ might not even have a physical meaning. The framework of two-scalable
functions guarantees that the iterations \eqref{eq: fixed point iterations}
converge to a fix point also in the case of type-II standard interference
functions, but the optimality meaning of this fixed point is no longer clear. Therefore,
in the following we consider a function $\ib$ of type-II and
assume it comes from a problem of the form \eqref{op: power control}.
With the framework of Fast-Lipschitz optimization
we characterize type-II standard power control problems to show that the fixed point
is also optimal for optimization in the form \eqref{op: power control}.
This is an important result that we can establish by Fast-Lipschitz
optimization.

As in the of standard functions in Section \eqref{sec: standard is FL},
we examine the problem in logarithmic variables $\x=\ln\p$ and arrive at the equivalent problem \eqref{op: type1 fl}, with
gradients given by \eqref{eq: gradients}.

\begin{theorem}
\label{thm: type-II-ish is FL} Assume $\ib(\p)$ be differentiable
and type-II standard, and consider $\fx=\ln\ib(e^{\x}).$ Let $\B=\left[B_{ij}\right]$
such that
\[
B_{ij}=\max_{\x}\left|\nabla_{i}f_{j}\left(\x\right)\right|=\max_{\p}\left|\nabla_{i}\i_{j}(\p)\frac{p_{i}}{\i_{j}\left(\p\right)}\right|,
\]
and assume $\rho\left(\B\right)<1$.  Let $\c>\0$ be an arbitrary
(positive) vector in $\Rn$ and let $\h(z)\in\Re^{m}$ be any strictly
increasing function of one variable. Then, problem~\eqref{op: power control}
is Fast-Lipschitz if
\begin{equation}\label{eq: s}
\s=\left(\I-\B\right)^{-1}\c
\end{equation}
 and
\begin{equation}\label{eq: pcc}
\pcc\left(\p\right)=\h\left({\textstyle \prod_{i}}p_{i}^{s_{i}}\right).
\end{equation}
\end{theorem}
\begin{IEEEproof}
Since $\rho\left(\B\right)<1,$ $\left(\I-\B\right)^{-1}$ is invertible
and
\[
\s=\left(\I-\B\right)^{-1}\c=\sum_{k=0}^{\infty}\B^{k}\c>\0.
\]
Let $\S=\text{diag}\left(\s\right)$ and introduce the scaled variables
$\y=\S\x$. The inverse of $\S$ exists and is positive since $\s>\0$.
Inserting $\x=\S^{-1}\y$ in~\eqref{op: type1 fl} gives the equivalent
problem
\begin{equation}
\OPmin{\y}{\go\left(\y\right)}{\y\ge\g\left(\y\right)=\S\,\f\left(\S^{-1}\y\right),}\label{op: type-II proof-1}
\end{equation}
where
\begin{align*}
\go\left(\y\right) & =\fo\left(\S^{-1}\y\right)=\h\left({\textstyle \prod_{i}}e^{y_{i}}\right)=\h\left(e^{\sum_{i}y_{i}}\right):=\h\left(z(\y)\right).
\end{align*}
The gradients of the problem are
\begin{align*}
\left[\nabla\go\left(\y\right)\right]_{ij} & =h_{j}'\big(z\left(\y\right)\big)\frac{\partial z\left(\y\right)}{\partial y_{i}}=h_{j}'\big(z\left(\y\right)\big)z\left(\y\right)
\end{align*}
and $\nabla\g(\y)=\S^{-1}\nabla\f(\x)\S$, where $\x=\S^{-1}\y$.

We will show that problem~\eqref{op: type-II proof-1} is Fast-Lipschitz
by \qc{Q2}. Condition \qcc{Q2.pos} is fulfilled because $\ib(\p)>\0$
and $\nabla\ib(\p)\le\0$ when $\ib$ is type-II standard. This gives
\[
\nabla_{i}g_{j}(\y)=\nabla_{i}f_{j}(\x)\frac{s_{j}}{s_{i}}=\left(\nabla_{i}\i_{j}(\p)\frac{p_{i}}{\i_{j}(\p)}\right)\frac{s_{j}}{s_{i}}\le0\quad\forall i,j,
\]
since $\s>\0$ and $\p=e^{\x}=e^{\S^{-1}\y}\ge\0$.

Condition\textbf{ }\qcc{Q2.inf}\textbf{ }requires\textbf{ }
\begin{equation}
\norm{\nabla\g\left(\y\right)}_{\infty}<\min_{j}\frac{\min_{i}\left[\nabla\go\left(\y\right)\right]_{ij}}{\max_{i}\left[\nabla\go\left(\y\right)\right]_{ij}},\label{eq: type-II temp 1-1}
\end{equation}
which is true by construction. To see this, note that
\[
\min_{j}\frac{\min_{i}\left[\nabla\go\left(\y\right)\right]_{ij}}{\max_{i}\left[\nabla\go\left(\y\right)\right]_{ij}}=\min_{j}\frac{\min_{i}h_{j}'\left(z\left(\y\right)\right)z\left(\y\right)}{\max_{i}h_{j}'\left(z\left(\y\right)\right)z\left(\y\right)}=1.
\]
The left side of \eqref{eq: type-II temp 1-1} therefore requires
\otc{%
\begin{align*}
\norm{\nabla\g\left(\y\right)}_{\infty} & =\norm{\S^{-1}\nabla\f(\x)\S}_{\infty}
=\max_{i}\sum_{j}\left|\nabla_{i}f_{j}\left(\x\right)\frac{s_{j}}{s_{i}}\right|<1
\end{align*}}{%
\begin{align*}
\norm{\nabla\g\left(\y\right)}_{\infty} & =\norm{\S^{-1}\nabla\f(\x)\S}_{\infty}
\\&
=\max_{i}\sum_{j}\left|\nabla_{i}f_{j}\left(\x\right)\frac{s_{j}}{s_{i}}\right|<1
\end{align*}}
for all $\y$ and $\x=\S^{-1}\y$. Since $s_{i}>0$ and $\left|\nabla_{i}f_{j}\left(\x\right)\right|\le B_{ij}$
for all $\x,$ this holds if
\[
\max_{i}\sum_{j}B_{ij}\frac{s_{j}}{s_{i}}<1,
\]
or equivalently, if $\sum_{j}B_{ij}s_{j}<s_{i}$ for all $i$. This
is the $i$th row of $\left(\I-\B\right)\s>\0$, which holds by construction
since
\[
\left(\I-\B\right)\s=\left(\I-\B\right)\left(\I-\B\right)^{-1}\c=\c>\0.
\]
Finally, condition \qcc{Q2.f0} is easily checked because
\[
\left[\nabla\go\left(\y\right)\right]_{ij}=h_{j}'\left(z\left(\y\right)\right)z\left(\y\right)>\0,
\]
which holds due to that $h_{j}\left(z\right)$ is
increasing and $z\left(\y\right)=e^{\sum_{i}y_{i}}>0$. This concludes
the proof.
\end{IEEEproof}

The form $\pcc\left(\p\right)=\h(z\left(\p\right))$, with $\h$ being
an increasing function, implies that all cost functions $\pcc$ that
can be handled with Theorem~\ref{thm: type-II-ish is FL} are equivalent
to the scalar cost $\kappa_{0}\left(\p\right)={\textstyle \prod_{i}}p_{i}^{s_{i}}$
obtained when $\h(z)=z$. If one instead chooses $\h\left(z\right)=\ln z$
the cost becomes $\kappa\left(\p\right)=\s^{T}\ln\p$, i.e., a weighted
sum of the power logarithms. Equation~\eqref{eq: s} states that
weighting $\s$ should lie in the interior of the cone spanned by
the columns of $\left(\I-\B^{T}\right)^{-1}.$

The assumption $\rho\left(\B\right)<1$ is crucial for Theorem~\ref{thm: type-II-ish is FL}
to hold, since it guarantees the existence of a positive scaling matrix
$\S$. The assumptions assure that $\rho\left(\nabla\fx\right)\le\norm{\nabla\fx}_{1}<1$,
so $\rho\left(\B\right)<1$ surely holds if there is a point $\x^{B}$
such that $\B=\left|\nabla\f\left(\x^{B}\right)\right|$. This means
that all elements of $\nabla\fx$ are minimized at the common point
$\x^{B}$. The simplest case where this is true is when $\nabla\fx=\A^{T}$
is constant. This requires an $\ib\left(\p\right)$ of the form
\otc{%
\begin{align*}
\i_{i}\left(\p\right) &
=\exp\left(\A\ln\p+\b\right)
=e^{b_{i}}\exp\left(\sum_{j}A_{ij}\ln p_{j}\right)
=e^{b_{i}}\exp\left(\sum_{j}\ln p_{j}^{A_{ij}}\right)
\\&
=e^{b_{i}}\exp\left(\ln\prod_{j}p_{j}^{A_{ij}}\right)
=e^{b_{i}}\prod_{j}p_{j}^{A_{ij}},
\end{align*}}{%
\begin{align*}
\i_{i}\left(\p\right) &
=\exp\left(\A\ln\p+\b\right)
=e^{b_{i}}\exp\left(\sum_{j}A_{ij}\ln p_{j}\right)
\\&
=e^{b_{i}}\exp\left(\sum_{j}\ln p_{j}^{A_{ij}}\right)
=e^{b_{i}}\exp\left(\ln\prod_{j}p_{j}^{A_{ij}}\right)
\\&
=e^{b_{i}}\prod_{j}p_{j}^{A_{ij}},
\end{align*}}
i.e., $\i_{i}\left(\p\right)$ should be a monomial. If problem~\eqref{op: power control}
has the basic cost function
\[
\kappa_{0}\left(\p\right)={\textstyle \prod_{i}}p_{i}^{s_{i}}
\]
from above, which also is a monomial, the problem is a geometric optimization
problem \cite{Boyd2007-2}. Interestingly, geometric problems become
convex with the change of variables $\x\triangleq\ln\p$, the same
variable transformation used throughout this section.

From the discussion above, we can establish the following new
qualifying condition:
\setcounter{qcn}{5}
\otc{\newlength\lblen
\newlength\qcondlen
\newlength\deflen
\newcommand{\moreroom}{\raisebox{-4pt}{\rule{0pt}{12pt}}}
\begingroup
\small
\setlength\lblen{30mm}
\setlength\qcondlen{110mm}
\begin{center}
\begin{tabular}{|c|@{}l|}

\hline
\multicolumn{2}{|c|}{\textbf{Qualifying Condition \ref{Qt2}}\moreroom}
\\ \hline 
\textbf{\qc{Qt2}} &
\begin{Tqcond}
  \label{Qt2}

  \case \label{Qt2.pos}
  $\nabla\f(\x)^2\ge\0, \quad \big($e.g., $\nabla\f(\x)\le\0\big)$

  \case \label{Qt2.inf}
   $ \big|\nabla\fx\!\big|\le\B$  and $\rho(\B)<1$

  \case \label{Qt2.f0}
  \begin{tabular}[t]{@{}l}
  $\fo(\x)= \h(\sum_i s_i x_i)$ for a strictly increasing function $\h(z)$\\
where $\s = (\I-\B)^{-1}\c$ and $\c>\0$
  \end{tabular}

\end{Tqcond}
\\
\hline
\end{tabular}
\end{center}
\endgroup }{\newlength\lblen
\newlength\qcondlen
\newlength\deflen
\newcommand{\moreroom}{\raisebox{-4pt}{\rule{0pt}{12pt}}}
\begingroup
\small
\setlength\lblen{30mm}
\setlength\qcondlen{72mm}
\begin{center}
\begin{tabular}{|c|@{}l|}

\hline
\multicolumn{2}{|c|}{\textbf{Qualifying Condition \ref{Qt2}}\moreroom}
\\ \hline 
\textbf{\qc{Qt2}} &
\begin{Tqcond}
  \label{Qt2}

  \case \label{Qt2.pos}
  $\nabla\f(\x)^2\ge\0, \quad \big($e.g., $\nabla\f(\x)\le\0\big)$

  \case \label{Qt2.inf}
   $ \big|\nabla\fx\!\big|\le\B$  and $\rho(\B)<1$

  \case \label{Qt2.f0}
  \begin{tabular}[t]{@{}l}
  $\fo(\x)= \h(\sum_i s_i x_i)$ for a strictly \\
  increasing function $\h(z)$ where \\
   $\s = (\I-\B)^{-1}\c$ and $\c>\0$
  \end{tabular}

\end{Tqcond}
\\
\hline
\end{tabular}
\end{center}
\endgroup }
\setcounter{qcn}{0}
Observe that the notation in \qcc{Qt2.inf} means the absolute value, not the norm. The new condition is numbered 6, although it is the 4th condition of this paper (see the appendix), but there are five known qualifying conditions in Fast-Lipschitz optimization~\cite{Jakobsson13}, the last two are not used in this paper and therefore are not included in the appendix.

\begin{theorem}
Assume problem \eqref{eq FLform} is feasible, and that qualifying
condition \qc{Qt2} above holds for every $\x\in\D$. Then, the problem
is Fast-Lipschitz, i.e., the unique Pareto optimal solution is given
by $\xa=\f(\xa).$ \end{theorem}
\begin{IEEEproof}
The proof is analogous to that of Theorem~\ref{thm: type-II-ish is FL}.
\end{IEEEproof}
The simplest form of the function $\fox$ in \qc{Qt2}\textbf{ }arises
from $\h(z)=z$. This means that $\fox=\s^{T}\x$ is a weighted sum
of the power logarithms, where the weights $\s$ are closely related
to the constraint gradient.

In the following section, we turn our attention to a class of power control problems that do not have monotonic constraint functions.

\section{Absolutely subhomogeneous interference functions \label{sec:Absolutely-subhomogeneous-interf}}
\newcommand\Ex{\mathbb{E}}
\newcommand\z{\mathbf{z}}
\newcommand\dd{\,\text{d}}
\newcommand\de{\text{d}}
\newcommand\zm{{z_{\min}}}
\newcommand\op{(\p)}

\newcommand\fad{\Theta}
\newcommand\pdf{\theta}
\newcommand\Nif{\Phi}

\newcommand\CM[1]{\textcolor[rgb]{0.44,0.00,0.94}{#1}}

In the previous sections we examined two-sided scalable functions
that where monotonically increasing (standard) and monotonically decreasing
(type-II standard). In the following we give an example of a problem
formulation where the constraints are not monotonic, hence neither
standard, nor type-II standard. We show convergence and optimality
through Fast-Lipschitz optimization, which was not known before.

The example builds upon the problem formulation in~\cite{Nuzman}.
Once again we consider problem~\eqref{op: power control} and assume that the cost function $\pcc\op$ is increasing in $\p$.
The formulation in \cite{Nuzman} starts with the affine SNR model~\eqref{linear interference function},
but adds a stochastic channel and outage as follows. Let
\begin{equation}
\label{nuz underlying IF}
\i_{i}(\p)=\frac{\sinr_{i}}{g_{ii}}\left(\sum_{j\neq i}g_{ij}p_{j}+\eta_{i}\right)
\end{equation}
 represent the \emph{expected} power needed to reach the SNR target
$\sinr_{i}$, and model the stochastic gain from transmitter $i$
to receiver $i$ by $g_{ii}\fad_{i}$ where $\fad_{i}$ is a stochastic
variable describing the fading of the wireless channel. Furthermore,
allow each transmitter to send only if the required power (to reach
the SNR target) is lower than some bound $b$. Combining the two effects
gives the new power control law
\begin{equation}
p_{i}^{k+1}=h\left(\frac{\i_{i}(\mathbf{p}^{k})}{\fad_{i}}\right),\label{eq: nuzman update 1}
\end{equation}
 where
\[
h(x)=\begin{cases}
x & \text{ if }x\le b,\\
0 & \text{ otherwise.}
\end{cases}
\]
 The fast timescale of the fading $\fad_{i}$ makes it hard to track
and measure in practice. Instead, let each transmitter node update its
transmit power according to the expectation \eqref{eq: nuzman update 1},
i.e.,
\begin{equation}\label{eq: nuz iter}
p_{i}^{k+1}=\Ex_{\fad_i}\left[h\left(\frac{\i_{i}(\mathbf{p}^{k})}{\fad_{i}}\right)\right]
\triangleq\Nif_{i}\left(\i_{i}\left(\p^{k}\right)\right).
\end{equation}
The expectation acts to smooth the discontinuous properties of $h(\cdot)$, and $\Nif_i(\i_i\op)$ is called the smoothed interference function of node (or mobile equipment) $i$.
The iterations in \eqref{eq: nuz iter} can be seen as a possible solution algorithm for a power control problem of the type
\begin{equation}
\label{smoth-upc}
 \begin{array}{ll}
 \min & \pcc(\p) \\
 {\rm s.t.}& p_{i}\ge f_i(\p) \triangleq \Nif_i(\i_i(\p)) \quad\forall i.
\end{array}
\end{equation}
However,  the nature of $h(x)$ will make $f_i(\p)$ non-monotonic, regardless
of underlying assumptions on $\i_i(\p)$. Therefore, neither the standard, nor the type-II standard interference
function approach applies here. To study the convergence properties of iterations based on these functions, \cite{Nuzman} introduces \emph{absolutely
subhomogeneous} functions, fulfilling
\[
e^{-\left|a\right|}\Nif(\x)\le\Nif(e^{a}\x)\le e^{\left|a\right|}\Nif(\x)
\]
for every $\x\ge\0$ and all scalars $a$. Note that absolute subhomogeneity
is implied by two-sided scalability.
In \cite{Nuzman} it is  shown that, if for each $i$,
 \begin{itemize}
   \item$ \i_i(\p)$ is \emph{standard}, and
   \item $\Nif_i(x)=\Ex_{\fad_i}\left[ h\left( x/\fad_i \right) \right]$ is bounded and absolutely subhomogeneous,
 \end{itemize}
then the sequence \eqref{eq: nuz iter} will converge to a fixed point.
However, nothing is said in \cite{Nuzman} about the
optimality of this fixed point.

Our approach is to use Fast-Lipschitz optimization
and qualifying condition \qc{Qinf}, which has no requirements on
the monotonicity of $\f\ox$.
Consider again problem~\eqref{smoth-upc}.
If
\[
\f(\p)=[f_1(\x),\,\dots,\,f_n(\x)]^T
\]
and $\pcc\left(\p\right)$ fulfill \qcc{Qinf.inf}, i.e.,
if $\nabla\pcc\left(\p\right)>\0$ and
\begin{equation*}
\norm{\nabla\f\left(\p\right)}_{\infty}<\frac{\qr\left(\p\right)}{1+\qr\left(\p\right)} \,,\label{eq: cond}
\end{equation*}
where
\[
\qr\left(\p\right)=\min_{j}\frac{\min_{i}\nabla_{i}\kappa_{j}\left(\p\right)}{\max_{i}\nabla_{i}\kappa{}_{j}\left(\p\right)},
\]
then problem~\eqref{smoth-upc} is Fast-Lipschitz and the iterations~\eqref{eq: nuz iter}
will converge to the optimal solution of \eqref{smoth-upc}. In the
previous sections, we used properties of standard and type-II standard
functions to show that the gradient norm $\normi{\nabla\f}$ was small
enough. In this section, we obtain the bound directly from $\normi{\nabla\ib}$
by using the following result:

\begin{lemma}\label{nuz lemma}
Let $\pdf_j(y)$ be the pdf of the channel
fading coefficient $\fad_j$, consider $z > 0$ and define
\begin{equation}
\label{nuz omega}
\Omega_j(z) \triangleq \int_{z}^{\infty}\frac{\pdf_j(y)}{y}dy-\pdf_j(z).
\end{equation}
Then, the infinity norm of the constraint function of problem~\eqref{smoth-upc} and the
infinity norm of the underlying interference function $\ib(\p)$ in~\eqref{nuz underlying IF} fulfill
\begin{equation}\label{eq: res}
\norm{\nabla\mathbf{f}(\p)}_{\infty}\le\max_{i}\left|\Omega_{i}\big(\i_{i}(\p)/b\big)\right|\,\norm{\nabla\i(\p)}_{\infty}.
\end{equation}
\end{lemma}

\begin{IEEEproof}
Dropping the explicit $\p$-dependence of $\i_j$ and $f_j$, we have
\otc{%
\begin{align*}
f_j & =\mathds{E}_{\fad_j}[h(\i_j/\fad_j)]=\int_{0}^{\infty} h(\i_j/y) \pdf_j(y) \dd y
\\ &
= \int_{0}^{\i_j/b} \underbrace{h(\i_j/y)}_{=0} \pdf_j(y) \dd y
+ \int_{\i_j/b}^{\infty} \underbrace{h(\i_j/y)}_{=\i/y} \pdf_j(y) \dd y
=\i_j \int_{\i_j/b}^{\infty} \frac{\pdf_j(y)}{y} \dd y,
\end{align*}}{%
\begin{align*}
f_j & =\mathds{E}_{\fad_j}[h(\i_j/\fad_j)]=\int_{0}^{\infty} h(\i_j/y) \pdf_j(y) \dd y
\\ &
= \int_{0}^{\i_j/b} \underbrace{h(\i_j/y)}_{=0} \pdf_j(y) \dd y
+ \int_{\i_j/b}^{\infty} \underbrace{h(\i_j/y)}_{=\i/y} \pdf_j(y) \dd y
\\&
=\i_j \int_{\i_j/b}^{\infty} \frac{\pdf_j(y)}{y} \dd y,
\end{align*}}
because
\[
h(\i_j/y)=\begin{cases}
\i_j/y & \text{ if }\i_j/y\le b \iff y \ge \i_j/b,\\
0 & \text{ otherwise.}
\end{cases}
\]
It follows that
\otc{%
\begin{align*}
\frac{\de f_j}{\de\i_j} &
=\frac{\de}{\de\i_j}\left(\i_j\int_{\i_j/b}^{\infty} \frac{\pdf(y)}{y} \dd y\right)
=\int_{\i_j/b}^{\infty} \frac{\pdf_j(y)}{y} \dd y - \pdf(\i_j/b)
\triangleq \Omega_j(\i_j/b).
\end{align*}}{%
\begin{align*}
\frac{\de f_j}{\de\i_j} &
=\frac{\de}{\de\i_j}\left(\i_j\int_{\i_j/b}^{\infty} \frac{\pdf(y)}{y} \dd y\right)
=\int_{\i_j/b}^{\infty} \frac{\pdf_j(y)}{y} \dd y - \pdf(\i_j/b)
\\ &
\triangleq \Omega_j(\i_j/b).
\end{align*}}
Returning to full notation, we have
\newcommand\ojp{\Omega_j\big(\i_j\op/b\big)}
\begin{align*}
\frac{\partial f_j\op}{\partial p_i}
&
= \frac{\de f_j\op}{\de \i_j\op}  \frac{\partial \i_j\op}{\partial p_i}
=  \ojp \frac{\partial \i_j\op}{\partial p_i}.
\end{align*}
It follows that
\otc{%
\begin{align*}
\normi{\nabla \f\op} &
= \max_i \sum_j\left|\frac{\partial f_j\op}{\partial p_i}\right|
= \max_i \sum_j\left| \ojp \frac{\partial \i_j \op}{\partial p_i} \right|
\\&
\le \max_j |\ojp| \cdot \max_i \sum_j\left| \frac{\partial \i_j(\p)}{\partial p_i} \right|
\\&
= \max_j |\ojp| \cdot \normi{\nabla\bm\i(\p)},
\end{align*}}{%
\begin{align*}
\normi{\nabla \f\op} &
= \max_i \sum_j\left|\frac{\partial f_j\op}{\partial p_i}\right|
\\&
= \max_i \sum_j\left| \ojp \frac{\partial \i_j \op}{\partial p_i} \right|
\\&
\le \max_j |\ojp| \cdot \max_i \sum_j\left| \frac{\partial \i_j(\p)}{\partial p_i} \right|
\\&
= \max_j |\ojp| \cdot \normi{\nabla\bm\i(\p)},
\end{align*}}
as is stated by~\eqref{eq: res}. This concludes the proof.
\end{IEEEproof}

Note that Lemma~\ref{nuz lemma} is true regardless of the underlying interference model $\ib(\p)$, e.g., $\ib\op$ does not need to be monotonic. We will use Lemma~\ref{nuz lemma} in a simplified form as follows:

\begin{corollary}
\label{cor: nuzman}
Suppose optimization problem \eqref{op: power control} fulfill qualifying condition \qcc{Qinf.inf} up to a scaling factor $\alpha>0$, i.e., if
\[
\alpha \normi{\i(\p)} < \frac{\qr(\p)}{1+\qr(\p)}.
\]
Then, optimization problem \eqref{smoth-upc} is Fast-Lipschitz if
\[
\max_{i,z} |\Omega_i(z)|  \le \alpha.
\]
\end{corollary}

This corollary allows us to say that problem~\eqref{smoth-upc}, regardless the underlying interference model $\ib\op$, is Fast-Lipschitz if
\begin{equation}
\label{nuz O req}
\max_{i,z} |\Omega_i(z)| <
\frac{1}{\normi{\i(\p)}}
\frac{\qr(\p)}{1+\qr(\p)}
\quad \forall\p.
\end{equation}
For fading coefficients from an arbitrary distribution, the function $\Omega_i(z)$ in equation~\eqref{nuz omega} might not be expressed on closed form. However, the max-value of $\Omega_i(z)$  can be found through numerical calculations.
We now apply Corollary~\ref{cor: nuzman} to two different distributions of the channel fading $\fad$, one is analyzed analytically and one is studied numerically.
\subsection{Fading models}
In what follows we consider two different fading models. First we investigate the case where the channel fading coefficient $\fad$ follows a Rayleigh distribution, whereby the worst-case value of $\Omega$ can be determined analytically. Thereafter, we investigate the case when $\fad$ follows an exponential distribution. In this case we find the worst-case value of  $\Omega$ through numeric calculation.
\subsubsection{Rayleigh distribution}
\newcommand\rdp{\lambda} 
Assume $\fad_i$ is follows a Rayleigh distribution with parameter $\rdp_i$ and with pdf
\begin{equation}
\label{nuz rayleigh}
\pdf_i(y)=
\frac{y}{\rdp_i^2}e^{-y^2/2\rdp_i^2},\quad\rdp_i>0.
\end{equation}
Recalling the definition of $\Omega_i(z)$ in \eqref{nuz omega}, we calculate the first term of $\Omega_i(z)$ as
\[
\int_{\z}^{\infty} \frac{\pdf_i(y)}{y} \dd y = \int_{\z}^{\infty} \frac{e^{-y^2/2\rdp_i^2}}{\rdp_i^2} \dd y.
\]
By the substitution $y=\sqrt2\rdp_i t$ we get $\dd y=\sqrt2\rdp_i \de t$ and
\begin{align*}
\int_{\z}^{\infty} \frac{\pdf_i(y)}{y} \dd y
&= \!
\int_{z/\sqrt2\rdp_i}^{\infty}  \!\!\!  \frac{ e^{-t^2}}{\rdp_i^2} \sqrt2\rdp_i \dd t
=\sqrt\frac{\pi}{2\rdp_i^2}\text{erfc}\left(\frac{z}{\sqrt2\rdp_i}\right),
\end{align*}
where $\text{erfc}(\cdot)$ is the complementary error function. Therefore, we have
\[
\Omega_i(z)=\sqrt\frac{\pi}{2\rdp_i^2}\text{erfc}\left(\frac{z}{\sqrt2\rdp_i}\right)
-\frac{z}{\rdp_i^2}e^{-z^2/2\rdp_i^2}.
\]
and
\[
\frac{\de \Omega_i(z)}{\de z}=\frac{e^{-z^2/2\rdp_i^2}}{\rdp_i^2}
\left(2-\frac{z^2}{\rdp_i^2}\right),
\]
which is smooth, and equal to zero only when $z=\sqrt2\rdp_i$. Therefore, the extreme values of $\Omega_i(z)$ must occur as $z\to0$, $\left.z=\sqrt2\rdp_i\right.$, or $z\to\infty$. Evaluating $\Omega_i$ at these points gives
\mbox{\(
\Omega_i(z)   \to   \sqrt{\pi/2}  \frac{1}{\rdp_i}
\)}
as
\(
z\to 0,
\)
\[
\Omega_i(\sqrt2\rdp_i)=\frac{1}{\rdp_i}
\underbrace{\left(\sqrt\frac{\pi}{2}\text{erfc}(1)-\sqrt2e^{-1} \right)}_{
\approx-0.323}>-\frac{1}{3\rdp_i},
\]
and
\(
\Omega_i(z)\to0\text{ as }z\to\infty
\)
 respectively. It follows that
 $\max_{i,z}|\Omega_i(z)|\le \alpha$
 if
\[
\alpha \ge \max\left \{\sqrt{\pi/2}  \frac{1}{\rdp_i} ,\, \frac13\frac1{\rdp_i} \right \}
\Leftrightarrow
\rdp_i \ge \frac{\sqrt{\pi/2}}{\alpha}
\]
for all $i$.
This means that if
\begin{itemize}
  \item[a)] the original (deterministic and outage-free) problem~\eqref{op: power control} is Fast-Lipschitz by qualifying condition \qc{Qinf}, i.e., $\alpha\le1$ in Corollary~\ref{cor: nuzman}, and
  \item[b)] the channel fading  $\fad_i$ follows a Rayleigh distribution~\eqref{nuz rayleigh} with parameter $\rdp_i\ge\sqrt{\pi/2}$,
\end{itemize}
then problem \eqref{smoth-upc} is Fast-Lipschitz by Corollary~\ref{cor: nuzman}. It follows that the iterations~\eqref{eq: nuz iter} converge to $\p^\star$, and $\p^\star$ is the unique optimal solution of the optimization problem~\eqref{smoth-upc}.

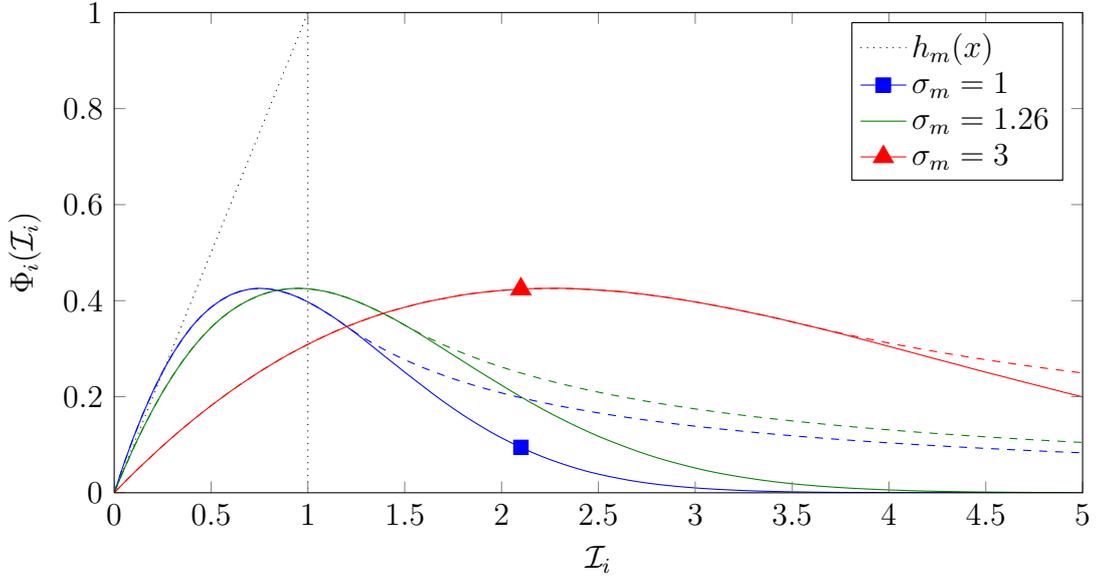
\begin{figure}
\centering
\iffastfig
        \framebox[\columnwidth][c]{\rule{0pt}{0.5\columnwidth}{Figure 1}}
\else
    \otc{
    \global\figwidth=0.8\columnwidth
    \global\figheight=0.5\figwidth
    }{
    \global\figwidth=0.8\columnwidth
    \global\figheight=0.5\figwidth
    }
    \tikzsetnextfilename{externalized_fig1}
     \begin{tikzpicture}
\definecolor{mycolor1}{rgb}{0,0.5,0}
\begin{axis}[%
scale only axis,
width=\figwidth,
height=\figheight,
xmin=0, xmax=5,
ymin=0, ymax=1,
xlabel={$\i_i$},
ylabel={$\Phi_i(\i_i)$},
axis on top,
legend entries={$h_m(x)$,$\sigma_m= 1$,$\sigma_m= 1.26$,$\sigma_m= 3$},
legend style={nodes=right}]

\addplot [
color=black,
dotted,
]
coordinates{
 (0,0)(1,1)(1,0)
};

\addplot [
color=blue,
mark=square*,
mark size = 0.5ex,
mark repeat=100,
mark phase=43,
solid
]
coordinates{
 (0,0)(0.05,0.0601259)(0.1,0.115499)(0.15,0.165708)(0.2,0.211141)(0.25,0.251867)(0.3,0.287836)(0.35,0.319044)(0.4,0.345914)(0.45,0.368235)(0.5,0.386331)(0.55,0.40075)(0.6,0.412049)(0.65,0.419518)(0.7,0.424096)(0.75,0.42571)(0.8,0.424454)(0.85,0.420811)(0.9,0.414923)(0.95,0.40715)(1,0.397416)(1.05,0.386265)(1.1,0.373927)(1.15,0.360746)(1.2,0.346356)(1.25,0.331322)(1.3,0.315517)(1.35,0.299417)(1.4,0.283129)(1.45,0.267253)(1.5,0.251075)(1.55,0.235559)(1.6,0.220071)(1.65,0.204639)(1.7,0.189808)(1.75,0.175708)(1.8,0.162113)(1.85,0.149203)(1.9,0.136991)(1.95,0.125087)(2,0.114182)(2.05,0.104117)(2.1,0.0945483)(2.15,0.0853744)(2.2,0.076868)(2.25,0.0691053)(2.3,0.0620489)(2.35,0.0554255)(2.4,0.0492243)(2.45,0.0436574)(2.5,0.0388422)(2.55,0.034449)(2.6,0.0303937)(2.65,0.0265859)(2.7,0.0232861)(2.75,0.0203765)(2.8,0.0177731)(2.85,0.0154768)(2.9,0.0134703)(2.95,0.0116909)(3,0.0101035)(3.05,0.00870293)(3.1,0.00756478)(3.15,0.00644414)(3.2,0.00550999)(3.25,0.00465642)(3.3,0.00390555)(3.35,0.0033139)(3.4,0.00284951)(3.45,0.00236751)(3.5,0.00197056)(3.55,0.00162797)(3.6,0.00136075)(3.65,0.00112587)(3.7,0.00093179)(3.75,0.000755157)(3.8,0.000642429)(3.85,0.000534076)(3.9,0.000426289)(3.95,0.000357239)(4,0.000285224)(4.05,0.000224335)(4.1,0.000166752)(4.15,0.000144623)(4.2,0.000114174)(4.25,8.73733e-05)(4.3,7.23039e-05)(4.35,6.3093e-05)(4.4,5.57478e-05)(4.45,4.42851e-05)(4.5,3.47288e-05)(4.55,2.70555e-05)(4.6,1.52863e-05)(4.65,1.3441e-05)(4.7,1.15678e-05)(4.75,1.16909e-05)(4.8,1.1814e-05)(4.85,5.89451e-06)(4.9,3.94131e-06)(4.95,3.98153e-06)(5,0)
};

\addplot [
color=mycolor1,
solid
]
coordinates{
 (0,0)(0.05,0.0481259)(0.1,0.0932151)(0.15,0.135136)(0.2,0.17395)(0.25,0.209749)(0.3,0.242448)(0.35,0.272364)(0.4,0.299168)(0.45,0.323159)(0.5,0.34427)(0.55,0.362585)(0.6,0.378219)(0.65,0.391361)(0.7,0.402092)(0.75,0.411092)(0.8,0.417569)(0.85,0.422135)(0.9,0.424874)(0.95,0.425791)(1,0.424961)(1.05,0.422183)(1.1,0.418512)(1.15,0.413154)(1.2,0.406689)(1.25,0.399005)(1.3,0.390423)(1.35,0.381098)(1.4,0.371181)(1.45,0.360553)(1.5,0.349113)(1.55,0.337461)(1.6,0.325061)(1.65,0.312288)(1.7,0.299638)(1.75,0.286671)(1.8,0.274133)(1.85,0.261164)(1.9,0.24859)(1.95,0.236294)(2,0.224004)(2.05,0.211595)(2.1,0.199743)(2.15,0.187942)(2.2,0.176796)(2.25,0.165908)(2.3,0.155454)(2.35,0.145423)(2.4,0.135824)(2.45,0.126514)(2.5,0.117543)(2.55,0.109453)(2.6,0.101549)(2.65,0.0939206)(2.7,0.0866341)(2.75,0.0797061)(2.8,0.0732667)(2.85,0.0673714)(2.9,0.0618516)(2.95,0.0565711)(3,0.0515191)(3.05,0.0468375)(3.1,0.0426187)(3.15,0.0388422)(3.2,0.0353016)(3.25,0.0319972)(3.3,0.0289021)(3.35,0.0260305)(3.4,0.0233923)(3.45,0.0210566)(3.5,0.0188473)(3.55,0.0169149)(3.6,0.0151511)(3.65,0.0135934)(3.7,0.0121521)(3.75,0.0108162)(3.8,0.00964957)(3.85,0.00853862)(3.9,0.00765509)(3.95,0.00673657)(4,0.00594437)(4.05,0.00523794)(4.1,0.00460607)(4.15,0.00399394)(4.2,0.00352478)(4.25,0.0030858)(4.3,0.00269984)(4.35,0.00234914)(4.4,0.00202003)(4.45,0.0017313)(4.5,0.00148741)(4.55,0.00131286)(4.6,0.00112011)(4.65,0.000965359)(4.7,0.000826948)(4.75,0.000707047)(4.8,0.000612005)(4.85,0.000535966)(4.9,0.000451033)(4.95,0.000379247)(5,0.000322812)
};

\addplot [
color=red,
solid,
mark=triangle*,
mark repeat=100,
mark phase=43,
mark size = .8ex
]
coordinates{
 (0,0)(0.05,0.0205953)(0.1,0.0406221)(0.15,0.0601259)(0.2,0.0790805)(0.25,0.0975993)(0.3,0.115499)(0.35,0.13275)(0.4,0.149423)(0.45,0.165708)(0.5,0.181318)(0.55,0.196349)(0.6,0.211141)(0.65,0.225203)(0.7,0.238675)(0.75,0.251867)(0.8,0.264259)(0.85,0.276376)(0.9,0.287836)(0.95,0.298658)(1,0.309023)(1.05,0.319044)(1.1,0.328482)(1.15,0.337642)(1.2,0.345914)(1.25,0.353752)(1.3,0.36119)(1.35,0.368235)(1.4,0.37492)(1.45,0.380821)(1.5,0.386331)(1.55,0.391595)(1.6,0.396535)(1.65,0.40075)(1.7,0.40455)(1.75,0.408367)(1.8,0.412049)(1.85,0.414833)(1.9,0.417346)(1.95,0.419518)(2,0.421414)(2.05,0.422795)(2.1,0.424096)(2.15,0.425018)(2.2,0.425412)(2.25,0.42571)(2.3,0.42574)(2.35,0.425353)(2.4,0.424454)(2.45,0.423483)(2.5,0.422183)(2.55,0.420811)(2.6,0.419245)(2.65,0.417317)(2.7,0.414923)(2.75,0.41253)(2.8,0.409808)(2.85,0.40715)(2.9,0.40395)(2.95,0.400836)(3,0.397416)(3.05,0.393776)(3.1,0.390122)(3.15,0.386265)(3.2,0.382252)(3.25,0.378197)(3.3,0.373927)(3.35,0.369877)(3.4,0.365368)(3.45,0.360746)(3.5,0.355878)(3.55,0.351036)(3.6,0.346356)(3.65,0.341427)(3.7,0.336461)(3.75,0.331322)(3.8,0.326078)(3.85,0.320836)(3.9,0.315517)(3.95,0.310183)(4,0.304805)(4.05,0.299417)(4.1,0.294044)(4.15,0.28856)(4.2,0.283129)(4.25,0.277823)(4.3,0.272564)(4.35,0.267253)(4.4,0.261733)(4.45,0.256398)(4.5,0.251075)(4.55,0.245929)(4.6,0.240875)(4.65,0.235559)(4.7,0.23024)(4.75,0.225178)(4.8,0.220071)(4.85,0.214847)(4.9,0.209606)(4.95,0.204639)(5,0.199743)
};

\addplot [
color=blue,
dashed
]
coordinates{
 (0,0)(0.05,0.0601259)(0.1,0.115499)(0.15,0.165708)(0.2,0.211141)(0.25,0.251867)(0.3,0.287836)(0.35,0.319044)(0.4,0.345914)(0.45,0.368235)(0.5,0.386331)(0.55,0.40075)(0.6,0.412049)(0.65,0.419518)(0.7,0.424096)(0.75,0.42571)(0.8,0.424454)(0.85,0.420811)(0.9,0.414923)(0.95,0.40715)(1,0.397416)(1.05,0.386265)(1.1,0.373927)(1.15,0.360746)(1.2,0.346356)(1.25,0.332502)(1.3,0.319713)(1.35,0.307872)(1.4,0.296877)(1.45,0.286639)(1.5,0.277085)(1.55,0.268147)(1.6,0.259767)(1.65,0.251895)(1.7,0.244487)(1.75,0.237501)(1.8,0.230904)(1.85,0.224663)(1.9,0.218751)(1.95,0.213142)(2,0.207814)(2.05,0.202745)(2.1,0.197918)(2.15,0.193315)(2.2,0.188921)(2.25,0.184723)(2.3,0.180707)(2.35,0.176863)(2.4,0.173178)(2.45,0.169644)(2.5,0.166251)(2.55,0.162991)(2.6,0.159857)(2.65,0.15684)(2.7,0.153936)(2.75,0.151137)(2.8,0.148438)(2.85,0.145834)(2.9,0.14332)(2.95,0.140891)(3,0.138542)(3.05,0.136271)(3.1,0.134073)(3.15,0.131945)(3.2,0.129884)(3.25,0.127885)(3.3,0.125948)(3.35,0.124068)(3.4,0.122243)(3.45,0.120472)(3.5,0.118751)(3.55,0.117078)(3.6,0.115452)(3.65,0.11387)(3.7,0.112332)(3.75,0.110834)(3.8,0.109376)(3.85,0.107955)(3.9,0.106571)(3.95,0.105222)(4,0.103907)(4.05,0.102624)(4.1,0.101372)(4.15,0.100151)(4.2,0.0989589)(4.25,0.0977946)(4.3,0.0966575)(4.35,0.0955465)(4.4,0.0944607)(4.45,0.0933994)(4.5,0.0923616)(4.55,0.0913466)(4.6,0.0903537)(4.65,0.0893822)(4.7,0.0884313)(4.75,0.0875005)(4.8,0.086589)(4.85,0.0856963)(4.9,0.0848219)(4.95,0.0839651)(5,0.0831254)
};

\addplot [
color=mycolor1,
dashed
]
coordinates{
 (0,0)(0.05,0.0481259)(0.1,0.0932151)(0.15,0.135136)(0.2,0.17395)(0.25,0.209749)(0.3,0.242448)(0.35,0.272364)(0.4,0.299168)(0.45,0.323159)(0.5,0.34427)(0.55,0.362585)(0.6,0.378219)(0.65,0.391361)(0.7,0.402092)(0.75,0.411092)(0.8,0.417569)(0.85,0.422135)(0.9,0.424874)(0.95,0.425791)(1,0.424961)(1.05,0.422183)(1.1,0.418512)(1.15,0.413154)(1.2,0.406689)(1.25,0.399005)(1.3,0.390423)(1.35,0.381098)(1.4,0.371181)(1.45,0.360553)(1.5,0.349113)(1.55,0.337852)(1.6,0.327294)(1.65,0.317376)(1.7,0.308041)(1.75,0.29924)(1.8,0.290928)(1.85,0.283065)(1.9,0.275616)(1.95,0.268549)(2,0.261835)(2.05,0.255449)(2.1,0.249367)(2.15,0.243567)(2.2,0.238032)(2.25,0.232742)(2.3,0.227683)(2.35,0.222838)(2.4,0.218196)(2.45,0.213743)(2.5,0.209468)(2.55,0.205361)(2.6,0.201412)(2.65,0.197611)(2.7,0.193952)(2.75,0.190425)(2.8,0.187025)(2.85,0.183744)(2.9,0.180576)(2.95,0.177515)(3,0.174557)(3.05,0.171695)(3.1,0.168926)(3.15,0.166244)(3.2,0.163647)(3.25,0.161129)(3.3,0.158688)(3.35,0.156319)(3.4,0.154021)(3.45,0.151788)(3.5,0.14962)(3.55,0.147513)(3.6,0.145464)(3.65,0.143471)(3.7,0.141532)(3.75,0.139645)(3.8,0.137808)(3.85,0.136018)(3.9,0.134274)(3.95,0.132575)(4,0.130917)(4.05,0.129301)(4.1,0.127724)(4.15,0.126186)(4.2,0.124683)(4.25,0.123216)(4.3,0.121784)(4.35,0.120384)(4.4,0.119016)(4.45,0.117679)(4.5,0.116371)(4.55,0.115092)(4.6,0.113841)(4.65,0.112617)(4.7,0.111419)(4.75,0.110246)(4.8,0.109098)(4.85,0.107973)(4.9,0.106871)(4.95,0.105792)(5,0.104734)
};

\addplot [
color=red,
dashed
]
coordinates{
 (0,0)(0.05,0.0205953)(0.1,0.0406221)(0.15,0.0601259)(0.2,0.0790805)(0.25,0.0975993)(0.3,0.115499)(0.35,0.13275)(0.4,0.149423)(0.45,0.165708)(0.5,0.181318)(0.55,0.196349)(0.6,0.211141)(0.65,0.225203)(0.7,0.238675)(0.75,0.251867)(0.8,0.264259)(0.85,0.276376)(0.9,0.287836)(0.95,0.298658)(1,0.309023)(1.05,0.319044)(1.1,0.328482)(1.15,0.337642)(1.2,0.345914)(1.25,0.353752)(1.3,0.36119)(1.35,0.368235)(1.4,0.37492)(1.45,0.380821)(1.5,0.386331)(1.55,0.391595)(1.6,0.396535)(1.65,0.40075)(1.7,0.40455)(1.75,0.408367)(1.8,0.412049)(1.85,0.414833)(1.9,0.417346)(1.95,0.419518)(2,0.421414)(2.05,0.422795)(2.1,0.424096)(2.15,0.425018)(2.2,0.425412)(2.25,0.42571)(2.3,0.42574)(2.35,0.425353)(2.4,0.424454)(2.45,0.423483)(2.5,0.422183)(2.55,0.420811)(2.6,0.419245)(2.65,0.417317)(2.7,0.414923)(2.75,0.41253)(2.8,0.409808)(2.85,0.40715)(2.9,0.40395)(2.95,0.400836)(3,0.397416)(3.05,0.393776)(3.1,0.390122)(3.15,0.386265)(3.2,0.382252)(3.25,0.378197)(3.3,0.373927)(3.35,0.369877)(3.4,0.365368)(3.45,0.360746)(3.5,0.355878)(3.55,0.351036)(3.6,0.346356)(3.65,0.341611)(3.7,0.336995)(3.75,0.332502)(3.8,0.328127)(3.85,0.323865)(3.9,0.319713)(3.95,0.315666)(4,0.31172)(4.05,0.307872)(4.1,0.304117)(4.15,0.300453)(4.2,0.296877)(4.25,0.293384)(4.3,0.289972)(4.35,0.286639)(4.4,0.283382)(4.45,0.280198)(4.5,0.277085)(4.55,0.27404)(4.6,0.271061)(4.65,0.268147)(4.7,0.265294)(4.75,0.262501)(4.8,0.259767)(4.85,0.257089)(4.9,0.254466)(4.95,0.251895)(5,0.249376)
};

\end{axis}
\end{tikzpicture}
\fi
\caption{The \emph{smoothed mobile behavior function} $\Nif(\i_i)$ for different $\rdp_i$, when $\fad_i$ follows a Rayleigh distribution. The dashed lines show the best approximations that are absolutely subhomogeneous, as required in \cite{Nuzman}. The dotted line shows the function $h(x)$ when $b=1$.}
\end{figure}
\hspace{1em}
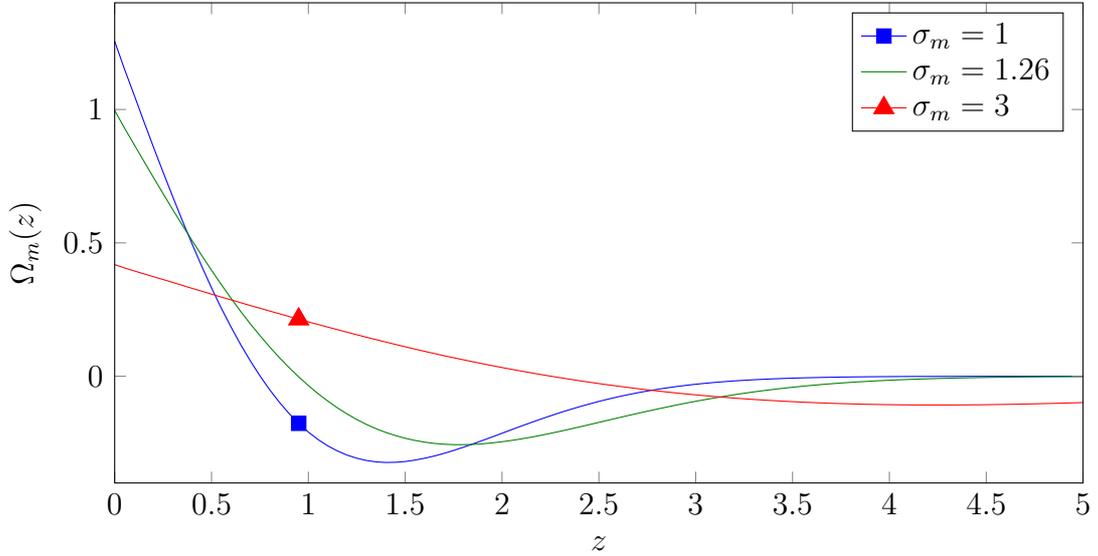
\begin{figure}
\centering
\iffastfig
        \framebox[\columnwidth][c]{\rule{0pt}{0.5\columnwidth}{Figure 2}}
\else
    \otc{
    \global\figwidth=0.8\columnwidth
    \global\figheight=0.5\figwidth
    }{
    \global\figwidth=0.8\columnwidth
    \global\figheight=0.5\figwidth
    }
    \tikzsetnextfilename{externalized_fig2}
%
%

\begin{tikzpicture}
\definecolor{mycolor1}{rgb}{0,0.5,0}

\begin{axis}[%
scale only axis,
width=\figwidth,
height=\figheight,
xmin=0, xmax=5,
ymin=-0.4, ymax=1.4,
ylabel={$\Omega_m(z)$},
xlabel={$z$},
axis on top,
legend entries={$\sigma_m= 1$,$\sigma_m= 1.26$,$\sigma_m= 3$},
legend style={nodes=right}]
\addplot [
color=blue,
solid,
mark=square*,
mark size = 0.5ex,
mark repeat=10000,
mark phase=20
]
coordinates{
 (0,1.25565)(0.05,1.15258)(0.1,1.05549)(0.15,0.956399)(0.2,0.859668)(0.25,0.76516)(0.3,0.672653)(0.35,0.582348)(0.4,0.495538)(0.45,0.411632)(0.5,0.331414)(0.55,0.255839)(0.6,0.185587)(0.65,0.119191)(0.7,0.0579582)(0.75,0.00148406)(0.8,-0.0503521)(0.85,-0.0972123)(0.9,-0.139253)(0.95,-0.176411)(1,-0.209114)(1.05,-0.237169)(1.1,-0.260749)(1.15,-0.279944)(1.2,-0.295473)(1.25,-0.307234)(1.3,-0.315719)(1.35,-0.320939)(1.4,-0.3232)(1.45,-0.322463)(1.5,-0.319595)(1.55,-0.314294)(1.6,-0.307315)(1.65,-0.298938)(1.7,-0.289117)(1.75,-0.278059)(1.8,-0.266155)(1.85,-0.253533)(1.9,-0.240401)(1.95,-0.227147)(2,-0.21358)(2.05,-0.199933)(2.1,-0.186503)(2.15,-0.173436)(2.2,-0.160688)(2.25,-0.148295)(2.3,-0.136335)(2.35,-0.124965)(2.4,-0.114213)(2.45,-0.104007)(2.5,-0.0943055)(2.55,-0.0852413)(2.6,-0.0768335)(2.65,-0.0690954)(2.7,-0.0619033)(2.75,-0.0552744)(2.8,-0.0492075)(2.85,-0.0436678)(2.9,-0.0386254)(2.95,-0.0340645)(3,-0.0299591)(3.05,-0.026273)(3.1,-0.0229447)(3.15,-0.0200174)(3.2,-0.0174014)(3.25,-0.015097)(3.3,-0.0130654)(3.35,-0.01126)(3.4,-0.00966354)(3.45,-0.00829268)(3.5,-0.0070932)(3.55,-0.00605215)(3.6,-0.00514373)(3.65,-0.00436189)(3.7,-0.0036878)(3.75,-0.00311297)(3.8,-0.00261179)(3.85,-0.00218828)(3.9,-0.00183272)(3.95,-0.00152599)(4,-0.00127054)(4.05,-0.00105557)(4.1,-0.000876686)(4.15,-0.000720643)(4.2,-0.000593359)(4.25,-0.000487796)(4.3,-0.000398537)(4.35,-0.000323966)(4.4,-0.000262425)(4.45,-0.000213047)(4.5,-0.000172576)(4.55,-0.000139439)(4.6,-0.000113606)(4.65,-9.0906e-05)(4.7,-7.25826e-05)(4.75,-5.74225e-05)(4.8,-4.52004e-05)(4.85,-3.66198e-05)(4.9,-2.91521e-05)(4.95,-2.28524e-05)(5,-1.86333e-05)
};

\addplot [
color=mycolor1,
solid
]
coordinates{
 (0,0.996546)(0.05,0.931049)(0.1,0.869361)(0.15,0.807093)(0.2,0.745353)(0.25,0.684593)(0.3,0.624478)(0.35,0.566068)(0.4,0.508349)(0.45,0.452197)(0.5,0.397446)(0.55,0.344292)(0.6,0.292944)(0.65,0.243681)(0.7,0.196553)(0.75,0.152407)(0.8,0.110043)(0.85,0.0701915)(0.9,0.0328319)(0.95,-0.00213964)(1,-0.0347463)(1.05,-0.0652807)(1.1,-0.0928503)(1.15,-0.118341)(1.2,-0.141358)(1.25,-0.162139)(1.3,-0.180569)(1.35,-0.196684)(1.4,-0.210539)(1.45,-0.222378)(1.5,-0.232417)(1.55,-0.240408)(1.6,-0.246855)(1.65,-0.251662)(1.7,-0.254687)(1.75,-0.25635)(1.8,-0.256375)(1.85,-0.255392)(1.9,-0.253089)(1.95,-0.249674)(2,-0.24542)(2.05,-0.240505)(2.1,-0.234715)(2.15,-0.228409)(2.2,-0.221411)(2.25,-0.214008)(2.3,-0.206213)(2.35,-0.198119)(2.4,-0.189803)(2.45,-0.181397)(2.5,-0.172942)(2.55,-0.164285)(2.6,-0.155756)(2.65,-0.147363)(2.7,-0.139119)(2.75,-0.131051)(2.8,-0.123142)(2.85,-0.115398)(2.9,-0.107901)(2.95,-0.100712)(3,-0.0938433)(3.05,-0.0872536)(3.1,-0.0809178)(3.15,-0.0748456)(3.2,-0.0691013)(3.25,-0.0636796)(3.3,-0.0585814)(3.35,-0.0537932)(3.4,-0.0493016)(3.45,-0.0450753)(3.5,-0.0411531)(3.55,-0.0374782)(3.6,-0.0340677)(3.65,-0.0308968)(3.7,-0.0279752)(3.75,-0.0252906)(3.8,-0.0228109)(3.85,-0.0205511)(3.9,-0.0184519)(3.95,-0.0165668)(4,-0.0148401)(4.05,-0.0132689)(4.1,-0.011843)(4.15,-0.0105632)(4.2,-0.00938807)(4.25,-0.0083336)(4.3,-0.00738374)(4.35,-0.00653266)(4.4,-0.00577406)(4.45,-0.00509487)(4.5,-0.00448602)(4.55,-0.00393468)(4.6,-0.0034532)(4.65,-0.00302275)(4.7,-0.00264212)(4.75,-0.00230541)(4.8,-0.00200631)(4.85,-0.00174158)(4.9,-0.00151282)(4.95,-0.00131169)(5,-0.00113439)
};

\addplot [
color=red,
solid,
mark=triangle*,
mark size = 0.8ex,
mark repeat=10000,
mark phase=20
]
coordinates{
 (0,0.418549)(0.05,0.406351)(0.1,0.395116)(0.15,0.384193)(0.2,0.373229)(0.25,0.362716)(0.3,0.35183)(0.35,0.340661)(0.4,0.329507)(0.45,0.3188)(0.5,0.307847)(0.55,0.296905)(0.6,0.286556)(0.65,0.275919)(0.7,0.265275)(0.75,0.255053)(0.8,0.24454)(0.85,0.234419)(0.9,0.224218)(0.95,0.213983)(1,0.203916)(1.05,0.194116)(1.1,0.184344)(1.15,0.174875)(1.2,0.165179)(1.25,0.15566)(1.3,0.146338)(1.35,0.137211)(1.4,0.128293)(1.45,0.119285)(1.5,0.110471)(1.55,0.101939)(1.6,0.0936248)(1.65,0.0852795)(1.7,0.0770997)(1.75,0.0693293)(1.8,0.0618622)(1.85,0.0542715)(1.9,0.0469088)(1.95,0.0397304)(2,0.0327653)(2.05,0.0258918)(2.1,0.0193194)(2.15,0.0128978)(2.2,0.00655773)(2.25,0.000494688)(2.3,-0.00537629)(2.35,-0.0111228)(2.4,-0.016784)(2.45,-0.0221783)(2.5,-0.0274179)(2.55,-0.0324041)(2.6,-0.037192)(2.65,-0.0418499)(2.7,-0.0464178)(2.75,-0.0507256)(2.8,-0.0548994)(2.85,-0.0588037)(2.9,-0.0626569)(2.95,-0.0662437)(3,-0.0697047)(3.05,-0.0730143)(3.1,-0.0761089)(3.15,-0.0790565)(3.2,-0.0818456)(3.25,-0.0844461)(3.3,-0.0869162)(3.35,-0.0891299)(3.4,-0.0912959)(3.45,-0.0933148)(3.5,-0.0952289)(3.55,-0.0969646)(3.6,-0.0984909)(3.65,-0.0999277)(3.7,-0.101221)(3.75,-0.102411)(3.8,-0.103486)(3.85,-0.10442)(3.9,-0.10524)(3.95,-0.105934)(4,-0.106515)(4.05,-0.10698)(4.1,-0.107326)(4.15,-0.107588)(4.2,-0.107733)(4.25,-0.107748)(4.3,-0.107657)(4.35,-0.107488)(4.4,-0.107279)(4.45,-0.106946)(4.5,-0.106532)(4.55,-0.106005)(4.6,-0.105389)(4.65,-0.104765)(4.7,-0.104079)(4.75,-0.103281)(4.8,-0.102438)(4.85,-0.10157)(4.9,-0.100658)(4.95,-0.0996459)(5,-0.0985804)
};

\end{axis}
\end{tikzpicture}
\fi
\vspace{-3mm}
\caption{This figure show the behaviour of $\Omega_i(z)$ for different $\rdp_i$. When $\alpha=1$, $\rdp_i=\sqrt{\pi/2}\approx1.26$ is the lower limit of $\rdp_i$ for which Corollary~\ref{cor: nuzman} applies (i.e., $|\Omega_i(z)|<1\;\forall z$).}
\label{figNuz2}
\vspace{0mm}
\end{figure}

\subsubsection{Exponential distribution}
\newcommand\edp{\lambda}

\begin{figure}
\centering
\iffastfig
        \framebox[\columnwidth][c]{\rule{0pt}{0.5\columnwidth}{Figure 3}}
\else
    \otc{
    \global\figwidth=0.8\columnwidth
    \global\figheight=0.5\figwidth
    }{
    \global\figwidth=0.8\columnwidth
    \global\figheight=0.5\figwidth
    }
    \tikzsetnextfilename{externalized_fig3}
    \input{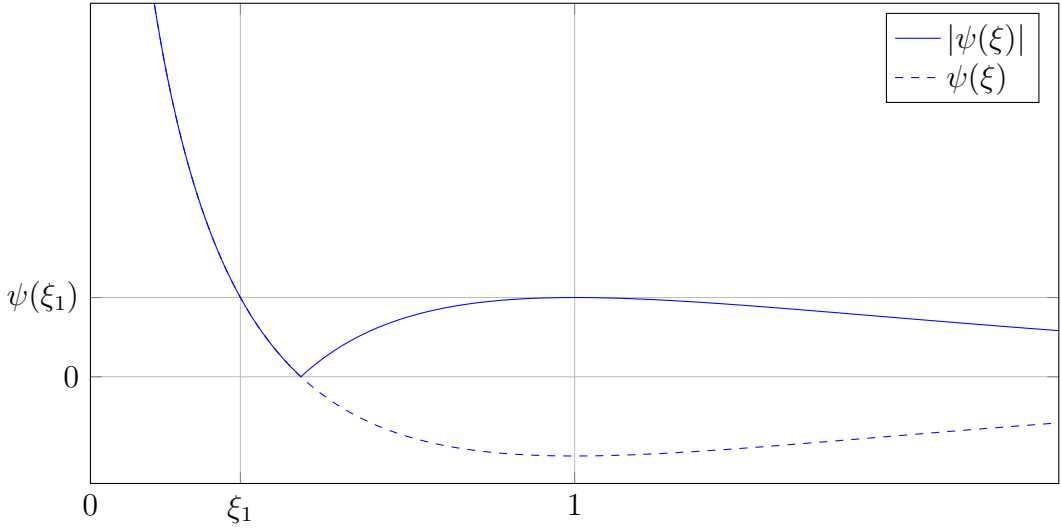}
\fi
\caption{Graph of $\psi(\xi)$ in equation \eqref{eq: psi}. }
\label{fig3}
\end{figure}

We now given an application of Corollary~\ref{cor: nuzman} to the case when the channel fading coefficients $\fad_{i}$
are exponentially distributed,
\begin{equation}
\label{nuz exponetial}
\fad_i\sim\pdf_i(y\,|\,\edp)=\edp e^{-\edp y}.
\end{equation}
This is know as Rayleigh fading. Denote $z\triangleq\i\op/b$ (we will drop the transmitter index $i$ to get a clearer notation), and highlight the $\edp$-dependence of $\Omega$ by writing
\otc{%
\begin{align*}
\Omega(z,\edp)&
= \int_{z}^{\infty}\frac{\pdf(y)}{y}\, dy-\pdf(z\,|\,\edp)
=\edp\!\left(\int_{\edp z}^{\infty}\frac{e^{-t}}{t}\, dt-e^{-\edp z}\!\right)
= \edp\psi(\edp z)
\end{align*}
}{%
\begin{align*}
\Omega(z,\edp)&
= \int_{z}^{\infty}\frac{\pdf(y)}{y}\, dy-\pdf(z\,|\,\edp)
=\edp\!\left(\int_{\edp z}^{\infty}\frac{e^{-t}}{t}\, dt-e^{-\edp z}\!\right)
\\&
= \edp\psi(\edp z)
\end{align*}}
where
\begin{equation}
\psi(\xi)\triangleq\int_{\xi}^{\infty}\frac{e^{-t}}{t}\, dt-e^{-\xi}.\label{eq: psi}
\end{equation}
The function $\psi(\xi)$ is shown in Figure~\ref{fig3}.

 To use the result in \eqref{nuz O req} we must show that the absolute
value of $\Omega$ is small enough. We will see that this is typically
the case, except when $z=\i/b$ goes to zero. This cannot happen in
practice, since the non-zero background noise $\eta$ always lower
bounds the interference. Therefore, we assume that $z=\i/b$ is lower
bounded by some $\zm$. For any given lower bound $z_{\min}$,
introduce
\[
\sigma_\zm\left(\edp\right)={\displaystyle \max_{z\ge z_{\min}}}\left|\Omega(z,\edp)\right|.
\]
The function $\sigma_{z_{\min}}$ is the worst case value over all possible values of $\edp$, given that $z\ge \zm$. To find $\sigma_{z_{\min}}\left(\edp\right)$, let
$\xi=\edp z$, whereby
\[
\left|\Omega(z,\edp)\right|
=\edp\left|\psi(\edp z)\right|
=\edp\left|\psi(\xi)\right|.
\]
For a fixed $\edp$, it is sufficient to find the
$z\ge z_{\min}$ that maximizes $\left|\psi(\edp z)\right|$ or,
equivalently, the $\xi\ge\xi_{\min}=\edp z_{\min}$ that maximizes
$|\psi(\xi)|$. Consider the plot of $\psi(\xi)$ is shown in Figure~\ref{fig3}. The derivative
\[
\frac{d\psi}{d\xi}=e^{-\xi}\left(1-\frac{1}{\xi}\right)
\]
 is zero only when $\xi=1$, and the second derivative is always positive.
The dashed lines highlight where $\xi=1$ and
\[
\xi=\xi_{1}=\left\{ t\,:\,\psi(t)=-\psi(1)\right\}.
\]
In order to maximize $|\psi(\xi)|$, it is clear that $\xi$ should be chosen as
\[
\xi = \begin{cases}
  1         & \text{if }    \xi_{1} \le \xi_{\min} \le 1, \\
  \xi_{\min} & \text{otherwise}.
\end{cases}
\]
In terms of the variables $\edp$ and $z$ we therefore have
\otc{%
\[
\sigma_{z_{\min}}\left(\edp\right)
={\displaystyle \max_{z\ge z_{\min}}}\left|\Omega(z,\edp)\right|
=\begin{cases}
\edp\psi(\edp z_{\min})=\Omega(z_{\min},\edp), & \text{if }\edp<\frac{\xi_{1}}{z_{\min}},\\
-\edp\psi(1), & \text{if }\frac{\xi_{1}}{z_{\min}}\le\edp\le\frac{1}{z_{\min}}\\
-\edp\psi(\edp z_{\min})=-\Omega(z_{\min},\edp), & \text{if }\frac{1}{z_{\min}}<\edp.
\end{cases},
\]}{%
\begin{multline}
\sigma_{z_{\min}}\left(\edp\right)
={\displaystyle \max_{z\ge z_{\min}}}\left|\Omega(z,\edp)\right|
\\
=\begin{cases}
\edp\psi(\edp z_{\min})=\Omega(z_{\min},\edp), & \text{if }\edp<\frac{\xi_{1}}{z_{\min}},\\
-\edp\psi(1), & \text{if }\frac{\xi_{1}}{z_{\min}}\le\edp\le\frac{1}{z_{\min}}\\
-\edp\psi(\edp z_{\min})=-\Omega(z_{\min},\edp), & \text{if }\frac{1}{z_{\min}}<\edp.
\end{cases},
\end{multline}}

It is clear that any stationary point of $\sigma_{z_{\min}}$ must
also be a stationary point of $\Omega(z_{\min},\edp)$, with derivative
\begin{align}
\nonumber \frac{d\Omega(z_{\min},\edp)}{d\edp} & =\frac{d}{d\edp}\left(\edp\psi\left(\edp z_{\min}\right)\right)=\psi\left(\edp z_{\min}\right)+\edp\frac{d\psi\!\left(\edp z_{\min}\right)}{d\edp}
\\& \nonumber
=\psi\left(\edp z_{\min}\right)+\edp\left(\zm e^{-\edp\zm}(1-\frac{1}{\edp\zm})\right)\\
&= \nonumber
\psi\left(\edp\zm\right)+e^{-\edp\zm}(\edp\zm-1).
\end{align}
 Setting the expression above to zero and solving numerically gives
the two solutions
\[
\begin{cases}
\edp\zm=v_{1}\approx0.1184 & \text{and}\\
\edp\zm=v_{2}\approx1.5656,
\end{cases}
\]
 i.e., when $\edp=v_{1}/\zm$ and $\edp=v_{2}/\zm$. Inserting
these values into $\sigma_{z_{\min}}\left(\edp\right)$ gives the
values
\[
\begin{cases}
\sigma_{z_{\min}}\left({\displaystyle \frac{v_{1}}{z_{\min}}}\right)={\displaystyle \frac{v_{1}}{z_{\min}}}\psi(v_{1})\approx{\displaystyle \frac{0.093}{z_{\min}}} & \text{and}\\
\rule{0pt}{20pt}\sigma_{z_{\min}}\left({\displaystyle \frac{v_{2}}{z_{\min}}}\right)={\displaystyle \frac{v_{2}}{z_{\min}}}\psi(v_{2})\approx{\displaystyle \frac{0.185}{z_{\min}}}
\end{cases}
\]
of the two local maxima shown in the Figure~\ref{fig4}.
\begin{figure}
\centering
\iffastfig
    \framebox[\columnwidth][c]{\rule{0pt}{0.5\columnwidth}{Figure 4}}
\else
    \otc{
    \global\figwidth=0.8\columnwidth
    \global\figheight=0.5\figwidth
    }{
    \global\figwidth=0.8\columnwidth
    \global\figheight=0.5\figwidth
    }
    \tikzsetnextfilename{externalized_fig4}
    \input{./tikz/fig4.tikz}
\fi
\caption{Plot of $\sigma_{z_{\min}}\left(\edp\right)$, note that the x-axis
is scaled by $\zm.$}
\label{fig4}
\end{figure}
Assuming $z\ge\zm$, we therefore have
\[
\left|\Omega(z,\edp)\right|
\le
\max_{z\ge \zm}\left|\Omega(\edp,z)\right|
=
\sigma_{z_{\min}}\left(\edp\right)
\le
0.185/\zm
\]
for any parameter value $\edp$ of the fading coefficient distribution parameter.
In particular, Corollary~\ref{cor: nuzman} states that problem~\eqref{smoth-upc} is Fast-Lipschitz if
\[
\norm{\nabla\ib\left(\p\right)}_{\infty} < \frac{0.185}{\zm} \frac{\qr\left(\p\right)}{1+\qr\left(\p\right)}
\]
for all $\p\ge\0$, where
\[
\qr\left(\p\right)=\min_{j}\frac{\min_{i}\nabla_{i}\kappa_{j}\left(\p\right)}{\max_{i}\nabla_{i}\kappa{}_{j}\left(\p\right)}
\]
is given by the characteristics of the cost function $\pcc\op$.

This example has showed how problems without monotonicity properties
can be analyzed with Fast-Lipschitz optimization. The price one has
to pay to ensure optimality is the tighter bound on $\norm{\nabla\ib\left(\p\right)}_{\infty}$
(note that $\qr\left(\p\right)/(1+\qr\left(\p\right))\le1/2$), as
opposed to requiring $\norm{\nabla\ib\left(\p\right)}_{\infty}<1$
for monotonic functions, which is sufficient to show contractivity.

\section{Conclusions and future work\label{sec:Conclusions-and-future}}

In this paper we examined the conditions under which power control
algorithms with standard, type-II standard and more general functions
fall under the Fast-Lipschitz framework. This allowed us to give the
studied problems a richer notion of optimality. In the process we
established a new qualifying condition for Fast-Lipschitz optimization
that shows a close relation between requirements on the cost function
and requirements on the constraints to achieve optimality.

In this paper we assumed that the functions are everywhere differentiable.
This is not necessarily required by the standard and type-II standard
formulations, and we believe this requirement can be dropped also
in Fast-Lipschitz optimization. However, this is something that still
needs to be formalized. Furthermore, the results of Section~\ref{sec: type-II is FL}
hint of possible relaxations of the qualifying conditions if one considers
cones different from the non-negative orthant.

\section*{Appendix: Fast-Lipschitz qualifying conditions}

Given a problem on Fast-Lipschitz form \eqref{eq FLform}, the General
Qualifying Condition (\qc{Qk}) of Table~\ref{tab:oldcond} guarantees that the problem is Fast-Lipschitz.
\begin{theorem}
[{\cite[Theorem 7]{Jakobsson13}}] \label{theorem: Main theorem - New}Assume
problem \eqref{eq FLform} is feasible, and that the General Qualifying
Conditions \qc{Qk} in Table~\ref{tab:oldcond} hold for every $\x\in\D$.
Then, the problem is Fast-Lipschitz, i.e.,
\[
\x^{k+1}:=\f(\x^{k})
\]
converges to $\xa=\f(\xa)$, and $\xa$ is the unique Pareto optimal
solution of problem \eqref{eq FLform}.
\end{theorem}
There are several special cases of \qc{Qk} that more convenient easier
to work with. We list three of them (the ones used in this paper) in
Table~\ref{tab:oldcond}.
\begin{proposition}
[\cite{Jakobsson13}]If any of qualifying conditions \qc{Q1}-\qc{Qinf}
hold, then so does \qc{Qk}.
\end{proposition}
\begin{remark}
Note that the qualifying conditions only are sufficient, not necessary.
This means that there can be problems that are Fast-Lipschitz but
fail to fulfill the qualifying conditions of Table~\ref{tab:oldcond}.
\end{remark}

\otc{
\newcommand\tableplace{h}
\otc{%
}{%
\def\arraystretch{1.5}
}
\begin{center}
\begin{table}[ht]
\renewcommand{\thetable}{\arabic{table}}
\centering
\small
\setcounter{qcn}{-1}
\setlength\lblen{25mm}
\setlength\qcondlen{0.9\linewidth}
\addtolength{\qcondlen}{-1\lblen}
\setlength\deflen{\qcondlen}
\begin{tabular}{|c|@{}l|}
\hline
\multicolumn{2}{|c|}{\textbf{General Qualifying Condition}}
\\ \hline
 \textbf{\qc{Qk}} &
\begin{Tqcond}
\itemsep=1.3\itemsep
  \label{Qk}
  \case $\nabla\fo(\x)\ge\0$ with non-zero rows
  \label{Qk.f0}
  \case
  \label{Qk.cont}
  $\norm{\nabla\f(\x)}<1$ 
  \case[There exists a $k \in \{1,2,\dots\}\cup\infty$ such that ]
  \label{Qk.pos}
 \begin{tabular}[t]{@{}l}
  when $k<\infty$, then \\
  $ \nabla\f(\x)^k\ge\0$
 \end{tabular}
  \case
  \label{Qk.inf}
  \begin{tabular}[t]{@{}l}
  when $k>1$, then \\
  $  \normi{\sum_{l=1}^{k-1} \nabla\fx^l} <  \q  \ox \triangleq \min_j \displaystyle\frac{\min_i [\nabla\fox]_{ij}}{\max_i [\nabla\fox]_{ij}}$
  \end{tabular}
\end{Tqcond}
 \\ \hline

\multicolumn{2}{|c|}{\textbf{Qualifying Condition \ref{Q1}}}
\\ \hline
 \textbf{\qc{Q1}} &
\begin{Tqcond}
\itemsep=1.3\itemsep
  \label{Q1}
  \case
  \label{Q1.f0}
  $\nabla\fo(\x)\ge\0$ with non-zero rows
  \case
  \label{Q1.cont}
  $ \norm{\nabla\fx}<1$
 \case
  \label{Q1.pos}
  $\nabla\f(\x)\ge\0$
\end{Tqcond}
\\ \hline

\multicolumn{2}{|c|}{\textbf{Qualifying Condition \ref{Q2}}}
\\ \hline
 \textbf{\qc{Q2}} &
 \begin{Tqcond}
\itemsep=1.3\itemsep
 \label{Q2}
  \case
  \label{Q2.f0}
  $\nabla\fo(\x)>\0$
  \case
  \label{Q2.pos}
  $\nabla\f(\x)^2\ge\0, \quad \big($e.g., $\nabla\f(\x)\le\0\big)$
  \case \label{Q2.inf}
  \begin{tabular}[t]{@{}l}
  \tabcolsep=0pt
  $\normi{ \nabla\fx} <  \q\ox \triangleq \min_j \displaystyle\frac{\min_i [\nabla\fox]_{ij}}{\max_i [\nabla\fox]_{ij}}$
  \end{tabular}
\end{Tqcond}
\\ \hline

\multicolumn{2}{|c|}{\textbf{Qualifying Condition \ref{Qinf}}}
\\ \hline
 \textbf{\qc{Qinf}} &
 \begin{Tqcond}
 \label{Qinf}
  \case
  \label{Qinf.f0}
  $\nabla\fo(\x)>\0$

  \case
  \label{Qinf.inf}
  $ \displaystyle \normi{\nabla\f(\x)}< \frac{\q\ox}{1+\q\ox}$
\end{Tqcond}
\\ \hline

\end{tabular}
\caption{Fast-Lipschitz qualifying conditions from~\cite{Jakobsson13}. Qualifying conditions~1-3 imply the general condition \qc{Qk}.}
\label{tab:oldcond}
\end{table}
\end{center}
\def\arraystretch{1} 
}{
\newcommand\tableplace{h}
\otc{%
}{%
\def\arraystretch{1.5}
}
\begin{center}
\begin{table}[h]
\renewcommand{\thetable}{\arabic{table}}
\centering
\small
\setcounter{qcn}{-1}
\setlength\lblen{25mm}
\setlength\qcondlen{1.05\linewidth}
\addtolength{\qcondlen}{-1\lblen}
\setlength\deflen{\qcondlen}
\begin{tabular}{|c|@{}l|}
\hline
\multicolumn{2}{|c|}{\textbf{General Qualifying Condition}}
\\ \hline
 \textbf{\qc{Qk}} &
\begin{Tqcond}
\itemsep=1.3\itemsep
  \label{Qk}
  \case $\nabla\fo(\x)\ge\0$ with non-zero rows
  \label{Qk.f0}
  \case
  \label{Qk.cont}
  $\norm{\nabla\f(\x)}<1$ 
  \case[There exists a $k \in \{1,2,\dots\}\cup\infty$ such that ]
  \label{Qk.pos}
 \begin{tabular}[t]{@{}l}
  when $k<\infty$, then \\
  $ \nabla\f(\x)^k\ge\0$
 \end{tabular}
  \case
  \label{Qk.inf}
  \begin{tabular}[t]{@{}l}
  when $k>1$, then \\
  $  \normi{\sum_{l=1}^{k-1} \nabla\fx^l} <  \q  (\x)$, \\
  where $\q\ox\triangleq \min_j \displaystyle\frac{\min_i [\nabla\fox]_{ij}}{\max_i [\nabla\fox]_{ij}}$
  \end{tabular}
\end{Tqcond}
 \\ \hline

\multicolumn{2}{|c|}{\textbf{Qualifying Condition \ref{Q1}}}
\\ \hline
 \textbf{\qc{Q1}} &
\begin{Tqcond}
\itemsep=1.3\itemsep
  \label{Q1}
  \case
  \label{Q1.f0}
  $\nabla\fo(\x)\ge\0$ with non-zero rows
  \case
  \label{Q1.cont}
  $ \norm{\nabla\fx}<1$
 \case
  \label{Q1.pos}
  $\nabla\f(\x)\ge\0$
\end{Tqcond}
\\ \hline

\multicolumn{2}{|c|}{\textbf{Qualifying Condition \ref{Q2}}}
\\ \hline
 \textbf{\qc{Q2}} &
 \begin{Tqcond}
\itemsep=1.3\itemsep
 \label{Q2}
  \case
  \label{Q2.f0}
  $\nabla\fo(\x)>\0$
  \case
  \label{Q2.pos}
  $\nabla\f(\x)^2\ge\0, \quad \big($e.g., $\nabla\f(\x)\le\0\big)$
  \case \label{Q2.inf}
  \begin{tabular}[t]{@{}l}
  \tabcolsep=0pt
  $\normi{ \nabla\fx} <  \q  (\x)$, \\
  where $\q\ox\triangleq \min_j \displaystyle\frac{\min_i [\nabla\fox]_{ij}}{\max_i [\nabla\fox]_{ij}}$
  \end{tabular}
\end{Tqcond}
\\ \hline

\multicolumn{2}{|c|}{\textbf{Qualifying Condition \ref{Qinf}}}
\\ \hline
 \textbf{\qc{Qinf}} &
 \begin{Tqcond}
 \label{Qinf}
  \case
  \label{Qinf.f0}
  $\nabla\fo(\x)>\0$

  \case
  \label{Qinf.inf}
  $ \displaystyle \normi{\nabla\f(\x)}< \frac{\q\ox}{1+\q\ox}$
\end{Tqcond}
\\ \hline

\end{tabular}
\caption{Fast-Lipschitz qualifying conditions from~\cite{Jakobsson13}. Qualifying conditions~1-3 imply the general condition \qc{Qk}, but are much easier to use from an analytical and computational point of view.}
\label{tab:oldcond}
\end{table}
\end{center}
\def\arraystretch{1} 
}

\bibliographystyle{IEEEtran}
\bibliography{references_ieee_formated,references}

\end{document}